\documentclass[12pt]{article} 
\usepackage{graphicx,amsmath, amssymb, verbatim}
\usepackage[all]{xy}
\newtheorem{proposition}{Proposition}

\newtheorem{definition}{Definition}

\newtheorem{conjecture}{Conjecture}

\newcommand {\C } {\mathbb{C}} 

\newcommand {\p } {\mathbb{P}} 
 
\newcommand {\Z} {\mathbb{Z}}

\newcommand{\ba}{\begin{eqnarray}}
\newcommand{\ea}{\end{eqnarray}}
\newcommand{\no}{\nonumber}

\begin{document} 

\title{\bf{J functions, non-nef toric varieties and equivariant local mirror symmetry of curves
}}

\author{Brian Forbes,\; Masao Jinzenji \\ \\ \it Division of 
Mathematics, Graduate School of Science \\ \it Hokkaido University \\ 
\it 
Kita-ku, Sapporo, 060-0810, Japan \\ {\it  brian@math.sci.hokudai.ac.jp }
{\it jin@math.sci.hokudai.ac.jp} }

\maketitle 

\begin {abstract} 
We develop techniques for computing the equivariant local mirror symmetry of curves, i.e. mirror symmetry for $\mathcal O(k)\oplus \mathcal O(-2-k)\rightarrow \p^1$ for $k \ge 1$. We also describe related methods for dealing with mirror symmetry of non-nef toric varieties. The basic tools are equivariant $I$ functions and their Birkhoff factorization.  
\end {abstract}

\section{Introduction}
In this paper, we discuss several problems related to the local mirror 
symmetry of $\mathcal O(k)$ $\oplus \mathcal O(-2-k)\rightarrow \p^{1}$, with particular emphasis on the cases where $k\geq 1$. Of course, this example has been discussed extensively 
by both mathematicians and physicists from the point of view of the 
topological vertex (e.g. \cite{BP}, \cite{OSV}), yet to date there has not been a discussion from the point of view of local mirror symmetry. 

The main difficulty in treating these cases via mirror symmetry comes from the fact that we often find ourselves working with non-nef toric manifolds (e.g. $F_n=\p(\mathcal O\oplus \mathcal O(-n))$ for $n\ge 3$). 
Recently, some of the techniques for carrying out the mirror computation of non-nef manifolds has been 
established (see for example \cite{CG}, \cite{I}, \cite{J3}). Using these results as a starting point, we develop methods for extracting Gromov-Witten invariants from mirror symmetry for the spaces in question.  

In particular, the example $\mathcal O(k)\oplus \mathcal O(-2-k)\rightarrow \p^1$ is also important in the geometrical interpretation of 
the integers obtained from genus 0 Gromov-Witten invariants of compact 
Calabi-Yau 3-folds $M$ via the multiple cover formula. These integers 
are heuristically considered as the number of rational curves in $M$, 
but since the normal bundle of a rational curve $C\subset M$ can be 
${\cal O}_{C}(k)\oplus{\cal O}_{C}(-2-k)$ for any $k \in \Z$, one 
expects that prepotential of 
${\cal O}(k)\oplus{\cal O}(-2-k)\rightarrow \p^{1}$ is the same as that of 
${\cal O}(-1)\oplus{\cal O}(-1)\rightarrow \p^{1}$. 
In this paper, we verify the following formula, which supports this 
naive speculation:
\begin{equation}
\int_{\bigl[\overline{M}_{0,0}(\p^1,d)\bigr]_{vir.}}
\biggl[\frac{c(R^{1}ft_{*}ev^{*}{\cal O}(-2-k),\lambda)}
{c(R^{0}ft_{*}ev^{*}{\cal O}(k),\lambda)}\biggr]_{2d-2}=
\frac{1}{d^{3}},\quad (k\geq 0),
\label{m1}
\end{equation} 
where $c(E,\lambda)$is the total Chern class 
$\sum_{j=0}^{rank(E)}\lambda^{j}c_{j}(E)$ and $[*]_{2d-2}$ represents the 
operation of picking out the coefficient of $\lambda^{2d-2}$. The l.h.s of 
(\ref{m1}) is derived from consideration of the degree $(d,0)$ genus 0 
Gromov-Witten invariant of $\p({\cal O}\oplus{\cal O}(k)\oplus
{\cal O}(-2-k))$. We now briefly outline this derivation.  
The moduli space of degree $(d,0)$ stable maps from $\p^{1}$ to 
 $\p({\cal O}\oplus{\cal O}(k)\oplus{\cal O}(-2-k))$ can be described as the 
projective bundle $\p({\cal O}\oplus R^{0}ft_{*}ev^{*}{\cal O}(k))$ 
over $\overline{M}_{0,0}(\p^{1},d)$. Therefore, the 0-point genus 0 Gromov-Witten 
invariant of degree $(d,0)$ is given by, 
\begin{equation}
\int_{\p({\cal O}\oplus R^{0}ft_{*}ev^{*}{\cal O}(k))}
c_{top}(R^{1}ft_{*}ev^{*}{\cal O}(-2-k)\otimes{\cal O}_{\p}(1)),
\label{top}
\end{equation}  
which directly leads us to the l.h.s. of (\ref{m1}) by the standard 
computation of Chern classes. In this paper, we also consider 
an integral which is the generalization of the l.h.s. of (\ref{m1}):
\begin{equation} 
\int_{\bigl[\overline{M}_{0,0}(\p^1,d)\bigr]_{vir.}}
\biggl[\frac{c(R^{1}ft_{*}ev^{*}{\cal O}(-2-k),\lambda)}
{c(R^{0}ft_{*}ev^{*}{\cal O}(k),z\lambda)}\biggr]_{2d-2},
\quad (k\geq 0).
\label{int2}
\end{equation}
If we set $z=-1$, the result turns out to be equivalent to 
the result of Bryan and Pandharipande \cite{BP} with anti-diagonal 
equivariant torus action. 

The main feature of our computation 
is the use of the $J$ function, obtained as the Birkhoff factorization of the $I$ function, which was invented by Coates and Givental
\cite{CG}. This approach turns out to be quite powerful for complex 3-dimensional manifolds.
Our basic method is to start from the twisted $I$ function of $\p^{1}$ , including the equivariant 
parameters that appear in Eqn. (\ref{int2}). We then take the asymptotical expansion of the $I$ function in the equivariant parameters, and Birkhoff factorize the result. 

What we find (for the diagonal torus action) is that in fact, the equivariant quantum cohomology of $\mathcal O(k)\oplus \mathcal O(-2-k)\rightarrow \p^1$ is the same for all $k\in \Z$:
\begin{conjecture}
\label{conjecture1}
Let $I_k^T(q)$ be the equivariant $I$ function of $X_k=\mathcal O(k)\oplus \mathcal O(-2-k)\rightarrow \p^1$, where $X_k$ is equipped with a $T^2$ action acting on the bundle $\mathcal O(k)\oplus \mathcal O(-2-k)$ with weights $(\lambda,\lambda)$. Then, for all $k\in \Z$, 
\begin{eqnarray}
	J_k^T(t)=e^{\lambda \tilde t_k/\hbar} I_{-1}^T(q),
\end{eqnarray}
where $J_k^T(t)$ is the $J$ function of $X_k$ after the shift by the mirror map, and $\tilde t_k$ is the `equivariant mirror map' read off from the coefficient of $\hbar^{-1}$ of the expansion of $J_k^T(t)$.
\end{conjecture}

We also discuss the quantum cohomology ring of $\p({\cal O}\oplus{\cal O}(k)\oplus
{\cal O}(-2-k))$ with motivation inspired by the result (\ref{m1}). Using similar techniques, we arrive at: 
\begin{conjecture}
\label{conjecture2}
 For all $k\in \Z$, 
$$QH^*(\p(\mathcal O \oplus \mathcal O(k) \oplus \mathcal O(-2-k))\cong QH^*(\p(\mathcal O \oplus \mathcal O(-1) \oplus \mathcal O(-1)).$$
\end{conjecture}  
We verify this conjecture by using the $J$-function together with the (ordinary) mirror map. We also present an alternative derivation, which relies on connection matrices on the moduli space.

The consideration of $\p(\mathcal O\oplus \mathcal O(k)\oplus \mathcal O(-2-k))$ directly leads us to another class of interesting examples of non-nef manifolds, namely the Hirzebruch 
surfaces $F_{n}\;\;(k\geq 3)$, which is the bundle $\p({\cal O}\oplus {\cal O}(-n))$ over $\p^{1}$. By again using the Birkhoff factorization, as well as the `generalized mirror transformation' (which means that we have a power series defining the mirror map for every generator of cohomology, rather than just the mirror map coming from 2 cycles), we arrive at the conjecture:
\begin{conjecture}
\label{conjecture3}
There are two isomorphism types of $QH^*(F_n)$, depending on whether $n$ is odd or even.
\end{conjecture} 
With this result, we then apply the technique of \cite{CG} to compute local mirror symmetry for $K_{F_3}$. 
 
The organization of this paper is the following. Section 2.1 is a short review of $I$ and $J$ functions for toric varieties, while 2.2-2.3 outline our techniques for describing (local) mirror symmetry for the spaces we are interested in. The remainder of the paper gives the details of applying the general theory to specific examples. Section 3 describes the equivariant local mirror symmetry of curves, while sections 4 and 5 deal with $\p(\mathcal O\oplus \mathcal O(k)\oplus \mathcal O(-2-k))$ and $F_n$, respectively.

\vspace{2em}
{\bf  Acknowledgment} We first would like to thank H. Iritani 
for helpful discussions and 
for generously giving us a computer program for Birkhoff factorization. We also 
thank Y. Konishi, M. van Manen, J. Bryan, M. Guest, S. Hosono, T. Eguchi and 
K. Liu for valuable discussions.  

The research of B.F. was 
funded by a COE grant of Hokkaido University. The research of M.J. is 
partially 
supported by JSPS grant No. 16740216.
\section{Overview}

In this section, we give a general guide to the computational strategies used throughout this paper. We begin with a brief review of Givental $I$ and $J$ functions, followed by our proposed methods of dealing with local mirror symmetry for curves and non-nef toric varieties, respectively. 

\subsection{Background}
\label{background}

Let $X$ be a compact K\"ahler toric variety with $\dim H_2(X)=k$. Note that we do not impose the condition $c_1(X)\ge 0$. Then $X$ can be described as a quotient $X=(\C^n-Z)/T^k$, where the weights of the torus action are given by an integral $k \times n$ matrix $M=(m_{ij})$, and $Z$ is the Stanley-Reisner ideal. We let $C_1\dots C_k$ be a basis of $H_2(X)$ corresponding to the rows of $M$, and choose K\"ahler classes $p_1\dots p_k \in H^{1,1}(X)$ satisfying $\int_{C_i} p_j=\delta_{ij}$. There are $n$ divisors $D_1\dots D_n$ in $X$ obtained by setting $z_i=0$, where $\C^n=( z_1  \dots  z_n )$. The intersection numbers between curves and divisors are then $C_i \cdot D_j =m_{ij}$. We denote the fundamental classes of these hypersurfaces by 
\begin{eqnarray}
    u_j=\sum_{i=1}^k m_{ij} p_i
\end{eqnarray}
which obey the relations $\int_C u_j=C \cdot D_j$ for any $C \in H_2(X)$. The first Chern class of $X$ is then $c_1(X)=\sum_{j=1}^n u_j$.

With these conventions, we have all the necessary ingredients for writing down the $I$ function
\begin{eqnarray}
	I_X=e^{(p_1\log q_1+\dots +p_k\log q_k)/\hbar}\sum_d \prod_{j=1}^n \frac{\prod_{m=-\infty}^{0}(u_j+m\hbar)}{\prod_{m=-\infty}^{\sum_i m_{ij} d_i}(u_j+m\hbar)}q_1^{d_1}\dots q_k^{d_k}.
\end{eqnarray}
 The coefficients take values in the cohomology ring of $X$:
\begin{eqnarray}
	H^*(X,\C)=\frac{\C[p_1 \dots p_k]}{\langle  u_{j_1}\dots u_{j_n} | (j_1\dots j_n)\subset Z  \rangle}.
\end{eqnarray}
The key feature of $I_X$ that we will make use of in this paper is that if $c_1(X)\ge 0$, then $I_X\in \C[[\hbar^{-1}]]$, but otherwise $I_X\in \C[[\hbar,\hbar^{-1}]]$. Now, in the $c_1(X)\ge 0$ case, we essentially have a structure theorem on the asymptotic form of $I_X$ \cite{CK} :
\begin{eqnarray}
\label{Iasymp}
	I_X = 1+\frac{\sum_{i=1}^k t_ip_i}{\hbar}+\frac{\sum_{\alpha \in H^4(X)} \alpha W_{\alpha}}{\hbar^2}+\dots
\end{eqnarray}
Above, $t_i$ give the mirror map, and the $W_{\alpha}$ are functions which can be used to compute Gromov-Witten invariants.

However, if $c_1(X)$ is not $\ge 0$, it is not at first glance clear how to derive the functions $t_i$, $W_{\alpha}$. Nonetheless, there is a closely related function, called the $J$ function, which does have the same nice structure.
We take the following as a definition of the $J$ function (from Corollary 5 of \cite{CG}):

\begin{definition}\label{birkhofffactor} \bf ($J$ function as Birkhoff factorization): \normalfont \it
Let $X$ be a compact K\"ahler toric manifold, and let $I_X(q,\hbar,\hbar^{-1}) \in \C[[\hbar,\hbar^{-1}]]$ be its corresponding $I$ function. Fix $c_{\alpha}(q,\hbar)\in \C[[\hbar]]$ such that 
\begin{eqnarray}
\label{Ibirkhoff}
\sum_{\alpha \in H^*(X)} c_{\alpha}(q,\hbar) \partial_{\alpha} I_X(q,\hbar,\hbar^{-1})=J_X(q,\hbar^{-1})\in \C[[\hbar^{-1}]].
\end{eqnarray}
Then we call $J_X(q,\hbar^{-1})$ a \normalfont $J$ function of $X$. 
\end{definition}
Notice that $J$ functions obtained as indicated will satisfy $J=I$ for any $X$ such that $c_1(X)\ge 0$. 

The main benefit in using $J_X$ is that often, we find that it has the same type of expansion (\ref{Iasymp}) as $I_X$ when $c_1(X)\ge 0$. 
Therefore, in nice cases, one would proceed by computing Gromov-Witten invariants and quantum cohomology by plugging the inverse mirror map into $J_X$. However, things are not always this simple, and we will explore possible complications in Section \ref{nonnef} below.

\subsection{Equivariant local mirror symmetry of curves}
\label{equivsection}

As our first application of Birkhoff factorization, we describe our proposed method for determining local mirror symmetry of curves. First, we need to be more precise about what we mean by `local mirror symmetry of curves'. 

Let $Y$ be a compact smooth K\"ahler toric variety of complex dimension 3, with an imbedded $\p^1 \hookrightarrow Y$. Then we can realize the normal bundle of $\p^1$ in $Y$ as a direct sum of line bundles, $N_{\p^1/Y}\cong \mathcal O(n_1)\oplus \mathcal O(n_2)=E_{n_1,n_2}$ for some $n_1,n_2\in \Z$. Then the question we are interested in is, what is the effective contribution of $\p^1 \hookrightarrow Y$ to the Gromov-Witten invariants of $Y$? 

We emphasize that this is \it not \normalfont the concept of mirror symmetry originally considered by Givental and Lian-Liu-Yau in \cite{G1}\cite{LLY}. The main difference is as follows. Let $\bar M_{0,0}(d,\p^1)$ denote the moduli stack of stable maps $f:\p^1\rightarrow \p^1$ with 0 marked points, such that $f_*[\p^1]=d\p^1$, $d\in \Z$. Define the usual forgetful and evaluation maps as $ft: \bar M_{0,1}(d,\p^1) \rightarrow \bar M_{0,0}(d,\p^1)$, $ev: \bar M_{0,1}(d,\p^1)\rightarrow \p^1$, obtained as forgetting and evaluation at the marked point, respectively. Then the constructions of \cite{G1}\cite{LLY} applied to the total space $E_{n_1,n_2}\rightarrow \p^1$ compute the moduli space integral 
\begin{eqnarray}
	\int_{[\bar M_{0,0}(d,\p^1)]_{vir}}c_{top}(R^0(ft_*ev^* E_{n_1,n_2}^+)\oplus R^1(ft_*ev^* E_{n_1,n_2}^- )),
\end{eqnarray}
which is a computation on the bundle $U_d \rightarrow \bar M_{0,0}(d,\p^1)$ with fiber $H^0(\p^1, f^* E_{n_1,n_2}^+)\oplus H^1(\p^1, f^* E_{n_1,n_2}^-)$. Here $E_{n_1,n_2}=E_{n_1,n_2}^+\oplus E_{n_1,n_2}^-$ is the splitting type of $E_{n_1,n_2}$, i.e. a separation into positive and negative bundles.

In contrast, the question we wish to address is the computation of the integral
\begin{eqnarray}
\label{modspaceint}
	\int_{[\bar M_{0,0}(d,\p^1)]_{vir}}\frac{c_{tot}(R^1(ft_*ev^* E_{n_1,n_2}^- ),\lambda)}{c_{tot}(R^0(ft_*ev^* E_{n_1,n_2}^+),z\lambda)}
\end{eqnarray}
where $c_{tot}(V,\lambda)=\sum_{j=0}^{rank(V)}\lambda^{rank(V)-j}c_{j}(V)$ denotes the total equivariant Chern class. The bundle in question here is the `direct difference' $H^1(\p^1, f^* E_{n_1,n_2}^-)\ominus H^0(\p^1, f^* E_{n_1,n_2}^+)$, which was defined rigorously in \cite{CG}.

In evaluating the integral (\ref{modspaceint}), we are interested in only one coefficient of the expansion of the integrand. Let $\lambda$ be an equivariant parameter. Then we have an expansion
\begin{eqnarray}
	\frac{c_{tot}(R^1(ft_*ev^* E_{n_1,n_2}^- ),\lambda)}{c_{tot}(R^0(ft_*ev^* E_{n_1,n_2}^+),z\lambda)}=\frac{\lambda^{d m_1-1}+\lambda^{dm_1-2}c_1(E_{n_1,n_2}^-)+\dots+c_{dm_1-1}(E_{n_1,n_2}^-)}{(z\lambda)^{dm_2}+(z\lambda)^{dm_2-1}c_1(E_{n_1,n_2}^+)+\dots+c_{dm_2}(E_{n_1,n_2}^+)}.
\end{eqnarray}
where $E_{n_1,n_2}^-,E_{n_1,n_2}^+$ have rank $m_1-1,m_2$ respectively.
Now, since $\bar M_{0,0}(d,\p^1)$ has dimension $2d-2$, we can show that the only term which is nonvanishing is the coefficient of $\lambda^1$ (this follows by setting $\lambda=1/\alpha$ and taking the coefficient of $\alpha^{2d-2}$). Therefore, it is this coefficient that we will be interested in computing via mirror symmetry.

With the above as background, we now detail our approach to local mirror symmetry on $X_{n_1,n_2}=E_{n_1,n_2}\rightarrow \p^1$. The starting point is the $I$ function for $\p^1$:
\begin{eqnarray}
	I_{\p^1}=e^{p\log q/\hbar}\sum_{d\ge 0}\frac{q^d}{\prod_{m=1}^d(p+m\hbar)^2}
\end{eqnarray}
with cohomology-valued coefficients (i.e. $p^2=0$). Then from \cite{CG}, the $I$ function of $X_{n_1,n_2}$ is obtained by twisting $I_{\p^1}$: 
\begin{eqnarray}
\label{twistedI}
	I_{n_1,n_2}=e^{p\log q/\hbar}\sum_{d\ge 0}\frac{q^2}{\prod_{m=1}^d(p+m\hbar)^2}\times \frac{\prod_{m=-\infty}^0(n_1p+m\hbar+\lambda_2)\prod_{m=-\infty}^0(n_2p+m\hbar+\lambda_1)}{\prod_{m=-\infty}^{n_1d}(n_1p+m\hbar+\lambda_2)\prod_{m=-\infty}^{n_2d}(n_2p+m\hbar+\lambda_1)},
\end{eqnarray}
where we set $\lambda_{1}=z\lambda,\;\;\lambda_{2}=\lambda$.
The main difficulty of this expression lies in using it to actually extract the relevant Gromov-Witten invariants. However, with the technology of Birkhoff factorization at our disposal, we can propose a means of overcoming this technicality.

As our first step, expand the series (\ref{twistedI}) in powers of $1/\lambda_1,1/\lambda_2$. For the sake of brevity, we set the equivariant parameters equal: $\lambda_1=\lambda_2=\lambda$. Unfortunately, the series expansion will introduce positive powers of $\hbar$ into the series for $I_{n_1,n_2}$, but as noted in Section \ref{background}, we can eliminate positive powers of $\hbar$ in the $I$ function by performing Birkhoff factorization. Let $J_{n_1,n_2}$ be the resulting Birkhoff factorized function. Then one can check directly that the series expansion for $J_{n_1,n_2}$ turns out to be 
\begin{eqnarray}
	1+\frac{pt(q)+\lambda \tilde t(q)}{\hbar}+\frac{\lambda p\big(W(q)+\tilde t(q) \log q\big)+\lambda^2 \tilde W(q)}{\hbar^2}+\dots
\end{eqnarray}
It is then straightforward to extract Gromov-Witten invariants from $J_{n_1,n_2}$. As usual, we interpret the coefficient of $\hbar^{-1}$ as the mirror map. Let $q(t)$ be the inverse of $t(q)$, and substitute this into $J_{n_1,n_2}$:
\begin{eqnarray}
	J_{n_1,n_2}(q(t))=1+\frac{pt+\lambda \tilde t(q(t))}{\hbar}+\frac{\lambda p\big(W(q(t))+
\tilde t(q(t))(\log(q(t))-t) +\tilde t(q(t)) t\big)+\lambda^2 \tilde W(q(t))}{\hbar^2}+\dots\no\\
\end{eqnarray}
We then still have a nontrivial component of the mirror map in the coefficient of $\hbar^{-1}$, the `equivariant mirror map', which we invert as:
\begin{eqnarray}
J'_{n_1,n_2}(t)&=&
	e^{-\lambda \tilde t(q(t))/\hbar}J_{n_1,n_2}(q(t))=1+\frac{pt}{\hbar}
\no\\
&&+\frac{\lambda p (W(q(t))
+\tilde t(q(t))(\log(q(t))-t))+\lambda^2\big(\tilde W(q(t))-\tilde t(q(t))^2/2\big)}{\hbar^2}+\dots
\end{eqnarray}
Then the function $J'_{n_1,n_2}(t)$ completely determines the equivariant quantum cohomology of $X_{n_1,n_2}$ (via Proposition \ref{jfn-qc} in the next subsection), and moreover we can directly extract Gromov-Witten invariants from the expansion of $W(q(t))$, the coefficient of $\lambda^1$, as expected. Moreover, in the applications we consider, the function $\tilde W(q(t))+\tilde t(q(t))^2/2$ turns out to have the same expansion as $W(q(t))$, up to an overall multiplicative constant.  

We will explore the application of this machinery to the Calabi-Yau case $n_1+n_2=-2$ later in the paper, and will find agreement with the recent results of Bryan-Pandharipande.

\subsection{Non-nef toric varieties}
\label{nonnef}

We now return to the discussion of the $J$ function as computed in Equation (\ref{Ibirkhoff}). We already mentioned that for well-behaved spaces, the asymptotic form of the $J$ function coincides with that of Eqn. (\ref{Iasymp}). However, in general, we may find the following phenomenon, which was first observed in \cite{J1}. Although the mirror map is indeed given as the coefficient of $1/\hbar$ in the expansion of $J_X$, the most general asymptotic expansion of the $J$ function is 
\begin{eqnarray}
	J_X=1+\hbar^{-1}\sum_{\alpha \in H^*(X)}\alpha t_{\alpha}+\dots
\end{eqnarray}
In nice cases, the sum above only runs over $\alpha \in H^2(X)$, but  this demonstrably fails for many spaces, notably $X=F_3=\p(\mathcal O\oplus \mathcal O(-3))$. Since $J_X$ is a function of $k$ variables $q_1\dots q_k$, and there are $N+1=\sum_j \dim H^{2j}(X)$ functions defining the mirror map, we are compelled to introduce a modified $J$ function $\hat J_X(q_0\dots q_N,\hbar^{-1})$, which possesses the extra variables we need to successfully invert the mirror map. 

A method of doing this was suggested in \cite{I}.  We will need to make the proposal of \cite{I} more concrete in order to carry out the computations we are interested in. The crucial ingredient are connection matrices, which will be constructed presently. Note that every object defined in this section, including the connection matrices, depends only on the $I$ function. 
\begin{definition}
Let $I_X$ be the $I$ function for a compact K\"ahler toric variety $X$, and let  $\{1,\alpha_1 \dots \alpha_{N}\}$ be a basis of $H^*(X)$. Define the fundamental solution
\begin{eqnarray}
	S^t=\begin{pmatrix}
	I_X & \partial_{\alpha_1} I_X & \dots & \partial_{\alpha_{N}} I_X
	\end{pmatrix}.
\end{eqnarray}
 Then the ($\hbar$ dependent) \normalfont connection matrices \it $\Omega_1\dots \Omega_k$ are defined by the equations 
\begin{eqnarray}
	\hbar q_i \partial/\partial q_i S= \Omega_i S, \ i=1\dots k.
\end{eqnarray}
\end{definition}
Unfortunately, these matrices $\Omega_i$ are not yet the ones which correspond to quantum multiplication by $p_i$ in the small quantum cohomology ring $QH^*(X)$. To compute the `right' connection matrices, we need their $\hbar$ independent form, which has been studied in \cite{Gu}\cite{I}. The first step is Birkhoff factorization of the fundamental solution: 
\begin{eqnarray}
	S(\hbar,\hbar^{-1})=Q(\hbar)R(\hbar^{-1}).
\end{eqnarray}
The positive part $Q(\hbar)$ then provides a gauge transformation which converts the $\Omega_i$ into $\hbar$ independent matrices:
\begin{eqnarray}
	\hat \Omega_i=Q^{-1}(\hbar)\Omega_i Q(\hbar)+\hbar q_i \partial/\partial q_i Q^{-1}(\hbar)Q(\hbar), \ i=1 \dots k.
\end{eqnarray}
Then the $\hat \Omega_1\dots \hat \Omega_k$ correspond to quantum multiplication by $p_1\dots p_k$ in $QH^*(X)$. We can immediately extend this to include quantum multiplication matrices for all $\alpha \in H^*(X)$; e.g., the operator corresponding to $p_1^2p_2$ is obviously $\hat \Omega_1^2 \hat \Omega_2$, etc. We let $\hat \Omega_{\alpha}$ denote the connection matrix corresponding to $\alpha \in H^*(X)$.

With these matrices in hand, we can write down the $J$ function for big quantum cohomology:
\begin{eqnarray}
	\hat J_X(q_0\dots q_{N},\hbar^{-1})=e^{\Theta/\hbar}J_X, \ \ \ \Theta=q_0+\sum_{j=k+1}^{N}\alpha_j q_j \hat \Omega_{\alpha_j}.
\end{eqnarray}
We remark that by using $\hat J_X(q_0\dots q_{N},\hbar^{-1})$ in the place of $I_X$ in Definition 2, connection matrices corresponding to multiplication in the big quantum cohomology ring can be computed.

With the modified $J$ function $\hat J_X(q_0\dots q_{N},\hbar^{-1})$ at our disposal, we are able to fully invert all functions $t_0\dots t_{N}$ of the mirror map. Let $q_0(t) \dots q_{N} (t)$ be the inverse mirror map. Then we find the coordinate shifted $J$ function in the following limit: 
\begin{eqnarray}
	J'_X(t_1\dots t_k,\hbar^{-1})=\lim_{t_0,t_{k+1}\dots t_{N}\rightarrow -\infty}\hat J_X(q_0(t)\dots q_{N}(t),\hbar^{-1})
\end{eqnarray}
Then $J'_X$ determines small quantum cohomology via the propostion:
\begin{proposition}
\label{jfn-qc}
Let $P(\hbar \partial/\partial t_i, e^{t_i},\hbar)$ be a differential operator such that $$P(\hbar \partial/\partial t_i, e^{t_i},\hbar)J'_X(t_1\dots t_k,\hbar^{-1})=0.$$ Then $P(p_i, e^{t_i},0)=0$ holds as a relation in small quantum cohomology.
\end{proposition}

Now that we have computed the correct $J$ function for $X$, $J'_X$, we can consider the effect of adding bundles to our case: $E\rightarrow X$. For simplicity, we here assume $E=\mathcal O(-\sum_{r=1}^k n_rp_r)$ with $n_i\ge 0 \ \forall i$. Expand 
\begin{eqnarray}
	J'_X=\sum_d J_d e^{dt}.
\end{eqnarray}
Then the \it twisted $J$ function \normalfont takes the form \cite{CG}:
\begin{eqnarray}
\label{Jtwist}
	J_E=\sum_d J_d\frac{\prod_{m=-\infty}^{0}(-\sum_{r=1}^k n_rp_r+m\hbar+\lambda)}{\prod_{m=-\infty}^{-\sum_{r=1}^k n_rd_r}(-\sum_{r=1}^k n_rp_r+m\hbar+\lambda)} e^{dt}
\end{eqnarray}
where $\lambda$ is an equivariant parameter. This function then gives Gromov-Witten invariants for the noncompact total space $E\rightarrow X$.

\section{Equivariant mirror symmetry for $\mathcal O(k)\oplus \mathcal O(-2-k)\rightarrow \p^1$}

We now turn our interest to local mirror symmetry for curves. The Gromov-Witten theory of curves has attracted attention recently due to its role in attractor equations \cite{OSV}. A recent paper of Bryan-Pandharipande has completely solved the $A$ model for all rank 2 bundles over a curve of arbitrary genus \cite{BP}. Here, we will see that at least some of their results can be reproduced easily from mirror symmetry, namely, the case where the base curve has genus 0. While we only actually solve the Calabi-Yau examples $\mathcal O(k)\oplus \mathcal O(-2-k)\rightarrow \p^1$, our method should work for arbitrary rank 2 bundles over $\p^1$, as explained in Section \ref{equivsection}. 

In this section, we will first present our evidence in favor of Conjecture \ref{conjecture1}. We then perform the computation for the antidiagonal action $\lambda_1=-\lambda_2$, giving the generating function of Gromov-Witten invariants for $\mathcal O(1)\oplus \mathcal O(-3)\rightarrow \p^1$.

\subsection{Equivariant Picard-Fuchs equations for $\mathcal O(-1)\oplus \mathcal O(-1)\rightarrow \p^1$}

Although the example $X_{-1}=\mathcal O(-1)\oplus \mathcal O(-1)\rightarrow \p^1$ has been studied extensively from a variety of perspectives, there has not yet been a satisfactory exposition which allows for generalization to all bundle spaces of type $X_k=\mathcal O(k)\oplus \mathcal O(-2-k)\rightarrow \p^1$. We will bridge this gap, and moreover give the expansion of the equivariant $I$ function $I_{-1}^T$ which will turn out to match $J_k^T$ (up to the equivariant mirror map) for all other $k$. 

Our discussion begins with the standard $I$ function for $X_{-1}$:
\begin{eqnarray}
	I_{-1}(q)=e^{p\log q/\hbar}\sum_{d\ge 0}\frac{\prod_{m=-d+1}^0(-p+m\hbar)^2}{\prod_{m=1}^d(p+m\hbar)^2}q^d.
\end{eqnarray}
As is well known, this series is annihilated by the differential operator 
\begin{eqnarray}
	\mathcal D_{-1}=\theta^2-q\theta^2, \ \theta= \hbar q\frac{d}{dq}
\end{eqnarray}
which has solution space $1, \log q$. This function is disappointingly free of instanton data, i.e. the multiple cover formula, if we consider it as a cohomology-valued hypergeometric series. 

We are interested in exhibiting the Gromov-Witten invariants of $X_{-1}$ in a way which generalizes to $X_k$ for all other $k$. The trick is to instead work with equivariant Gromov-Witten invariants, where we consider the effect of a torus action $(\lambda,\lambda)$ on the bundle $\mathcal O(-1)\oplus \mathcal O(-1)$. This yields the equivariant $I$ function 
\begin{eqnarray}
	I^T_{-1}=e^{p\log q/\hbar}\sum_{d\ge 0}\frac{\prod_{m=-d+1}^0(-p+m\hbar+\lambda)^2}{\prod_{m=1}^d(p+m\hbar)^2}q^d
\end{eqnarray}
which is annihilated by the \it equivariant differential operator \normalfont
\begin{eqnarray}
\label{equivDE1}
	\mathcal D_{-1}^T=\theta^2-q(\theta-\lambda)^2.
\end{eqnarray}
The interesting fact is that while $\mathcal D_{-1}f=0$ did not yield any instanton information, the equivariant equation $\mathcal D_{-1}^T f=0$ does indeed. This is seen most easily by directly expanding $I^T_{-1}$:
\begin{eqnarray}
\label{basicI0}
I^T_{-1}=1+\frac{p \log q}{\hbar}+\frac{-2 \lambda p Li_2(q) +\lambda^2Li_2(q)}{\hbar^2}+\frac{p\lambda^2 S_1 +\lambda^3 S_2}{\hbar^3}+\dots	
\end{eqnarray}
where $S_1,S_2$ are power series in $q$ whose exact form is not relevant here (although these functions do precisely match those of $J_k^T$ as indicated in Conjecture \ref{conjecture1}), and the polylogarithm function is 
\begin{eqnarray}
	Li_k(x)=\sum_{n>0}\frac{x^n}{n^k}.
\end{eqnarray}

We would like to point out that the differential operator of Equation (\ref{equivDE1}) should be thought of as a `remedy' for the insufficient Picard-Fuchs differential operator $\theta^2-q\theta^2$ which comes from local mirror symmetry. While the authors have constructed an extended system which overcomes this difficulty in a previous paper \cite{FJ1}, the above $\mathcal D_{-1}^T$ has the advantage that the classical limit $q\rightarrow 0$ reproduces the ordinary cohomology relation for $\p^1$. 

In summary, then, the two important aspects of this subsection are (1) we can use equivariant techniques to recover Gromov-Witten invariants which are not visible from the original geometry $X_{-1}$, and (2) the expression (\ref{basicI}), which we will use to compare with the result on other $X_k$.

\subsection{Nonvanishing invariants for $\mathcal O\oplus \mathcal O(-2)\rightarrow \p^1$}

We move on to the next most basic case, $X_0=\mathcal O\oplus \mathcal O(-2)\rightarrow \p^1$. It turns out that while Birkhoff factorization is not necessary, we will again need to introduce an equivariant parameter to the $I$ function in order to exhibit nonzero invariants. 

First, we recall `usual' nonequivariant mirror symmetry for $X_0$. The $I$ function for $X_0$ is given as:
\begin{eqnarray}
	I_0=e^{p \log q/\hbar}\sum_{d\ge 0} \frac{\prod_{m=-2d+1}^0(-2p+m\hbar)}{\prod_{m=1}^d(p+m\hbar)^2}q^d.
\end{eqnarray}
For the moment, we ignore the fact that the coefficients of $I_0$ are cohomology-valued, and attempt to compute Gromov-Witten invariants for $X_0$ by using the coefficient of $1/\hbar^2$ of $I_0$. The asymptotic expansion of $I_0$ is 
\begin{eqnarray}
	I_0=1+\frac{pt}{\hbar}+\frac{p^2W}{\hbar^2}+\dots.
\end{eqnarray}
Here $t$ is the usual mirror map 
\begin{eqnarray}
	t=\log q + 2\sum_{n\ge 1}\frac{(2n-1)!}{(n!)^2}q^n.
\end{eqnarray}
Then inserting the inverse mirror map $q(t)$ into $W$, we find $W=t^2p^2/2$, which means that all Gromov-Witten invariants are 0. There are two ways to see why this has to be true. First, $X_0$ is essentially a local $K3$ surface, which is known to have vanishing Gromov-Witten invariants. Secondly, closer examination of $I_0$ reveals that what is being computed is the moduli space integral 
\begin{eqnarray}
\label{modint}
	\int_{\bar M_{0,0}(\p^1,d)}c_{2d-1}(R^1 ft_* ev^*\mathcal O(-2))
\end{eqnarray}
where $\bar M_{0,0}(\p^1,d)$ is the moduli space of degree $d$ maps $f:\p^1\rightarrow \p^1, f_*[\p^1]=d[\p^1]$ with 0 marked points, and $ft$, $ev$ are the usual evaluation and forgetful maps $ev:\bar M_{0,1}(\p^1,d)\rightarrow \p^1$, $ft: \bar M_{0,1}(\p^1,d) \rightarrow \bar M_{0,0}(\p^1,d)$. Then since the dimension of $\bar M_{0,0}(\p^1,d)$ is $2d-2$, we must have that the integral (\ref{modint}) is zero. Moreover, we expect to be able to recover nonzero invariants by instead integrating $c_{2d-2}(R^1 ft_* ev^*\mathcal O(-2))$.

If one thinks in terms of equivariant mirror symmetry for curves, a natural way of proceeding becomes clear. We use the equivariant vector bundle $E^T\rightarrow \p^1$, obtained by considering a $T^2$ action on the bundle $\mathcal O \oplus \mathcal O(-2)$ with weights $(0,\lambda)$. Set $\mathcal U_d=R^1 ft_* ev^* E^T$. Then we can write the total equivariant Chern polynomial of $\mathcal U_d$ as 
\begin{eqnarray}
	c_{tot}(\mathcal U_d)=c_{2d-1}(\mathcal U_d)+\lambda c_{2d-2}(\mathcal U_d)+\dots +\lambda^{2d-1}.
\end{eqnarray}
This means that in principle, we ought to be able to find the `real' invariants of $X_0$ by examining the coefficient of $\lambda$ of the equivariant $I$ function  
\begin{eqnarray}
	I_0^T=e^{p \log q/\hbar}\sum_{d\ge 0} \frac{\prod_{m=-2d+1}^0(-2p+m\hbar+\lambda)}{\prod_{m=1}^d(p+m\hbar)^2}q^d.
\end{eqnarray}
It turns out that this function does in fact compute the expected invariants. 
We first restrict the coefficients to the cohomology ring $\C[p]/\langle p^2 \rangle$. Then expanding in powers of $1/\hbar$:
\begin{eqnarray}
	I_0^T=1+\frac{p t-\lambda \tilde t}{\hbar}+\frac{p\lambda(W-\tilde t\log q) +\lambda^2 \tilde W}{\hbar^2}+\dots
\end{eqnarray}
where $\tilde t=(t-\log q)/2$.
Then inserting the inverse mirror map into $I_0^T$, we have 
\begin{eqnarray}
J_0^T(t)=	1+\frac{p t+\lambda Li_1(e^t)}{\hbar}+\frac{-\lambda p (2Li_2(e^t)+t Li_1(e^t))+\lambda^2 \tilde W}{\hbar^2}+\dots.
\end{eqnarray}
For the last step, we have to invert the `equivariant mirror map', namely $\lambda Li_1(e^t)$. After doing this, we find 
\begin{eqnarray}
e^{-\lambda Li_1(e^t)/\hbar}J_0^T(t)=1+\frac{p t}{\hbar}+\frac{-2 p\lambda Li_2(e^t)+\lambda^2 Li_2(e^t)}{\hbar^2}+\dots	
\end{eqnarray}
which agrees exactly with Equation (\ref{basicI0}) for $X_{-1}$! We also checked higher powers of $1/\hbar$, which exhibit complete agreement between the two expressions.

From the perspective of differential operators, our calculations here imply the following. In considerations of local mirror symmetry, we should really be replacing the ordinary Picard-Fuchs operator for $X_0$, which annihilates $I_0$
\begin{eqnarray}
	\mathcal D_0=\theta^2-q 2\theta(2\theta+\hbar), \ \theta=\hbar q\frac{d}{dq}
\end{eqnarray}
with the annihilator of $I^T_0$, namely the equivariant Picard-Fuchs operator 
\begin{eqnarray}
	\mathcal D_0^T=\theta^2-q(2\theta-\lambda)(2\theta-\lambda+\hbar).
\end{eqnarray}

From the considerations of \cite{FJ1}, the above means that in principle, we should be able to reconstruct all of local mirror symmetry simply by examining equivariant $I$ functions. This would give a simple method of verifying physical Gromov-Witten invariants without concern for noncompactness of the toric variety.

\subsection{ $k\ge 1$}

We now generalize the above approach to include all $X_k=\mathcal O(k)\oplus \mathcal O(-2-k)\rightarrow \p^1$, $k\ge 1$. It is indeed possible to derive a version of mirror symmetry, in the sense that we have a mirror map and a double logarithmic function which reproduces the Gromov-Witten invariants of \cite{BP}. Comparing to the $X_0$ case, it turns out that $k\ge 1$ forces us to use the Birkhoff factorization to correctly compute invariants. We also find that the result is a collection of rational functions in the equivariant weights, so that we can only read off enumerative information by specifying values for the weights. 

We consider equivariant Gromov-Witten theory on $X_k$, endowed with a $T^2$ action with weights $(\lambda_1,\lambda_2)$ on the respective bundle factors $\mathcal O(k)\oplus \mathcal O(-2-k)$. Then \cite{CG} tells us that mirror symmetry for $X_k$ should be encoded in the following hypergeometric function:
\begin{eqnarray}
	I_k^T=e^{p\log q/\hbar}\sum_{d\ge 0}\frac{\prod_{m=(-2-k)d+1}^0((-2-k)p+m\hbar+\lambda_2)}{\prod_{m=1}^d(p+m\hbar)^2 \prod_{m=1}^{kd}(kp+m\hbar+\lambda_1)}q^d=e^{p\log q/\hbar}\sum_{d\ge 0}C_d^k(\lambda_1,\lambda_2)q^d.
\end{eqnarray}
While this may be the correct $I$ function, in its given form it is not clear exactly how one should extract Gromov-Witten data from it. For example, if $k=1$, then in the nonequivariant limit $\lambda_1=\lambda_2=0$, this function reduces to $I_{K_{\p^2}}$, the $I$ function for $\mathcal O(-3)\rightarrow \p^2$, whose invariants we are not presently interested in. Our guiding principle at the moment is that we would like to reproduce the famous multiple cover formula for curves which is known from calculations on $X_{-1}$.

There is a way around these difficulties, as we described in Section \ref{equivsection}, which goes as follows. We expand the coefficients $C_d^k(\lambda_1,\lambda_2)$ in powers of $1/\lambda_1$, which introduces positive powers of $\hbar$, and then Birkhoff factorize the resulting expression. We denote the result of Birkhoff factorization by $J_k^T$. 
Then the surprising fact is that $J_k^T$ actually contains exactly the expected mirror symmetry data! In other words, we find (up to order 2 in  $1/\hbar$): 
\begin{eqnarray}\no
	J_k^T(q)=1+\frac{p t^k(q,\lambda_1,\lambda_2)+\tilde t^k(q,\lambda_1,\lambda_2)}{\hbar}+\frac{p \big(W^k(q,\lambda_1,\lambda_2)+\tilde t^k(q,\lambda_1,\lambda_2)\log(q) \big)+\tilde W^k(q,\lambda_1,\lambda_2)}{\hbar^2}
\end{eqnarray}
Let $q(t)$ be the inverse of $t^k(q,\lambda_1,\lambda_2)$. As in the previous section, we first apply the inverse mirror map: 
\begin{eqnarray}\no
	J_k^T(t)=1+\frac{p t+ \tilde t^k(q(t),\lambda_1,\lambda_2)}{\hbar}+
\frac{p\big(\hat W^{k}(q(t),\tilde t^k(q(t)),\lambda_1,\lambda_2)+
\tilde t^{k}(q(t),\lambda_{1},\lambda_{2})t\big)+
\tilde W^k(q(t),\lambda_1,\lambda_2)}{\hbar^2}
\end{eqnarray}
where 
\begin{equation}
\hat W^k(q(t),\tilde t^k(q(t)),\lambda_1,\lambda_2):= 
W^k(q(t),\lambda_1,\lambda_2)+\tilde t^k(q(t))(\log q(t)-t),
\end{equation}
and then invert the equivariant mirror map $\tilde t^k(q(t),\lambda_1,\lambda_2)$:
\begin{eqnarray*}
	e^{-\tilde t^k(q(t),\lambda_1,\lambda_2)/\hbar}J_k^T(t)=1+\frac{pt}{\hbar}+\frac{p\hat W^k(q(t),\tilde t^{k}(q(t)),\lambda_1,\lambda_2)+\tilde W^k(q(t),\lambda_1,\lambda_2)-\tilde t^k(q(t),\lambda_1,\lambda_2)^2/2
}{\hbar^2}
\end{eqnarray*}
Then, if we specialize the torus weights so that $\lambda_1=\lambda_2=\lambda$, we find that for any $k$, 
\begin{eqnarray}
	e^{-\tilde t^k(q(t),\lambda)/\hbar}J_k^T(t)=1+\frac{p t}{\hbar}+\frac{-2 p\lambda Li_2(e^t)+\lambda^2 Li_2(e^t)}{\hbar^2}+\dots
\end{eqnarray}
which is exactly the formula we found for $I_{-1}^T$. This establishes our Conjecture, and agrees with the results of \cite{BP} for the diagonal torus action $\lambda_1=\lambda_2$.

As an example, we list here the resulting mirror map, equivariant mirror map and double log function for $X_1=\mathcal O(1)\oplus \mathcal O(-3)\rightarrow \p^1$ for the diagonal and antidiagonal torus action. The functions $W(q(t))$ are used to enumerate Gromov-Witten invariants.

\bigskip
\noindent
Diagonal action $\lambda_1=\lambda_2=\lambda$:
\begin{eqnarray}\no
t(q)&=&\log q +20q+536q^2+\frac{73280}{3}q^3+1404096q^4+92091392 q^5+\dots 
\\ \no
 \tilde t(q,\lambda)&=&-\lambda\big(4q+88q^2+\frac{10816}{3}q^3+193728q^4+\frac{60621824}{5}q^5\dots\big)
\\ \no
	W(q,\lambda)&=&-\lambda \big(2q+\frac{241}{2}q^2+\frac{48566}{9}q^3+\frac{6981379}{24}q^4+\frac{1344390356}{75}q^5\dots \big) \\
\hat W(q(t),\tilde t(q(t),\lambda),\lambda)&=&-2\lambda \Big(e^t+\frac{1}{4}e^{2t}+\frac{1}{9}e^{3t}+\frac{1}{16}e^{4t}+\frac{1}{25}e^{5t}+\dots\Big),
\end{eqnarray}

\bigskip
\noindent
Antidiagonal action $\lambda_1=-\lambda_2=\lambda$ (equivariantly Calabi-Yau case):
\begin{eqnarray}\no
	t(q)&=&\log q -8q+74q^2-\frac{3212}{3}q^3+18609q^4-\frac{1787308}{5}q^5+\dots \\ \no
	\tilde t(q,\lambda)&=&\lambda \big(2q-17q^2+\frac{710}{3}q^3-\frac{8049}{2}q^4+\frac{381142}{5}q^5\dots \big)
	\\ \no
	W(q,\lambda)&=& \lambda \big(4q-55q^2+\frac{7600}{9}q^3-\frac{179005}{12}q^4+\frac{21600262}{75}q^5-\dots \big)\\
\hat W(q(t),\tilde t(q(t),\lambda),\lambda)
&=&\lambda \big( 4e^t-7e^{2t}+\frac{220}{9}e^{3t}-
\frac{455}{4}e^{4t}+\frac{15504}{25}e^{5t}-\dots \big)
\end{eqnarray}
This concludes our discussion on the local mirror symmetry for curves.

We remark that in the process of actually carrying out this computation, the series expansion in $1/\lambda_1$ requires substantial computer power. As such, we have found difficulty in doing the calculation to much higher order. The advantage, however, is that we are able to make direct contact with mirror symmetry and the hypergeometric series $I_k^T$. At any rate, we have succeeded in our goal of describing mirror symmetry for $\mathcal O(k)\oplus \mathcal O(-2-k)\rightarrow \p^1$.

Next, we describe an alternate computational method which allows us to check the results to much higher order. Again, the starting point is an equivariant hypergeometric series.

\subsection{$k\ge 1$: the $A$ model computation}
In this subsection, we compute local Gromov-Witten invariants of 
${\cal O}(k)\oplus {\cal O}(-k-2)\rightarrow \p^{1}$ from the A-model  
point of view by using the fixed point theorem. 
We will see that the results 
completely agree with the ones obtained above from local mirror symmetry. 
This provides more evidence for the validity of the results found through our mirror computation. In the process of doing the 
fixed point computation, we will see appearance of ``dynamics'' in the
combinatorial structure of the mirror map. 

First, we define the total Chern class of a vector bundle $E$:   
\begin{equation}
c(E,\lambda):=\sum_{n=0}^{rank(E)}c_{n}(E)\lambda^{n},
\end{equation}
and introduce our definition of the generating function of local Gromov-Witten invariants of 
${\cal O}(k)\oplus {\cal O}(-k-2)\rightarrow \p^{1}$ as follows:
\begin{equation}
F(q,z)=\sum_{d=1}^{\infty}q^{d}
\int_{\bigl[\overline{M}_{0,0}(\p^1,d)\bigr]_{vir.}}
\biggl[\frac{c(R^{1}ft_{*}ev^{*}({\cal O}_{\p}(-2-k)),\lambda)}
{c(R^{0}ft_{*}ev^{*}({\cal O}_{\p}(k)),z\lambda)}\biggr]_{2d-2},
\label{adef}
\end{equation}
where $[*]_{2d-2}$ means the operation of picking out the coefficient 
of $\lambda^{2d-2}$.
In this setting, $z=1$ (resp. $z$=-1) corresponds to the diagonal 
(resp. antidiagonal) action in the previous computation. $F_{q,z}$ 
can be computed by using the Atiyah-Bott fixed-point theorem under 
the following $\C^{\times}$ action on $\p^1$:
\begin{equation}
T(X_{1}: X_{2}):=(\exp(u_{1}t)X_{1}: \exp(u_{2}t)X_{2}),
\label{torus}
\end{equation}
where we set $u_{1}=0,u_{2}=1$ in the following computation. 
By applying standard results (for example \cite{K}, \cite{J1}, \cite{CK}),
we can express the coefficient $q^{d}$ of $F(q,z)$ in terms of a summation 
of contributions of colored tree graphs with degree. Let $\Gamma$ be a 
colored graph with degree. Each edge $\alpha\in Edge(\Gamma)$ has degree 
$d_{\alpha}$ which takes on a positive integer value, 
and each vertex $v\in Vert(\Gamma)$
has color $c(v)\in\{1,2\}$ which corresponds to the two fixed points $(1:0)$ and 
$(0:1)$ of $\p^{1}$ under the $\C^{\times }$ action (\ref{torus}). Two colors 
$c(v_{1})$ and $c(v_{2})$ must be different if $v_{1}$ and $v_{2}$ are 
directly connected by a edge of $\Gamma$. We denote by $d_{v}$ the sum of 
degree $d_{\alpha}$  of $\alpha\in Edge(\Gamma)$ which is connected to $v$. 
Of course, we can define the degree of $\Gamma$ by 
$\sum_{\alpha\in Edge(\Gamma)}d_{\alpha}$, and denote by $G_{d}$ the set 
of all the colored tree graphs with degree $d$.
With this setup, we can write down the coefficient of $q^{d}$ in $F(q,z)$ 
in terms of the graph sum as follows:     
\begin{eqnarray}
&&\int_{\bigl[\overline{M}_{0,0}(\p^1,d)\bigr]_{vir.}}
\biggl[\frac{c(R^{1}ft_{*}ev^{*}({\cal O}_{\p}(-2-k)),\lambda)}
{c(R^{0}ft_{*}ev^{*}({\cal O}_{\p}(k)),z\lambda)}\biggr]_{2d-2}\nonumber\\
&&=\biggl[\sum_{\Gamma\in G_{d}}(-1)^{e(\Gamma)+d}
\frac{1}{|Aut(\Gamma)|}\prod_{\alpha\in Edge(\Gamma)}
d_{\alpha}\cdot a_{d_{\alpha}}\prod_{v\in Vert(\Gamma)}
(d_{v})^{val(v)-3}(f_{v})^{val(v)-1}\biggr]_{2d-2},
\label{gs}
\end{eqnarray}
where 
\begin{equation}
a_{d}:=\frac{d^{2d}}{(d!)^{2}}
\cdot\frac{\prod_{m=-(2+k)d+1}^{-1}(1+\frac{m}{d}\lambda)}
{\prod_{m=1}^{kd}(1+\frac{m}{d}z\lambda)},
\end{equation}
and 
\begin{eqnarray}
f_{v}:=\left\{\begin{array}{cc}1&
(\mbox{if}\;\;\; c(v)=1),\\
(1-(k+2)\lambda)(1+kz\lambda)
\;\;&(\mbox{if}\;\;\; c(v)=2). \end{array}\right.
\label{fv}
\end{eqnarray}
 For brevity, we introduce a polynomial $f$ that 
appears in (\ref{fv}) as follows:  
\begin{equation}
f:=(1-(k+2)\lambda)(1+kz\lambda).
\end{equation}
Therefore, we can compute $F(q,z)$ by summing up the contributions of 
all the colored tree graphs with degree. One problem is that the total number of graphs 
in $G_{d}$ increases dramatically as the degree $d$ rises, so we give 
the following technical discussion to address the question of efficiently summing up tree graphs.

First, we sum up the contributions from the graphs that have only one 
edge, which appear only once in $G_{d}$ for each $d$. Looking back at 
(\ref{gs}), we can sum up the contribution of these graphs as follows:
\begin{equation}
A(a_1,a_2,a_3,\cdots):=\sum_{d=1}^{\infty}\frac{(-1)^{d-1}a_{d}}{d^{3}}q^d.
\label{a}
\end{equation}
In (\ref{a}), the notation $A(a_1,a_2,a_3,\cdots)$ means that we regard 
$A$ as function of the $a_{d}$'s. Next, we consider star graphs. By a 
`star graph', we mean any graph whose vertices are connected to only one edge (i.e.,
$val(v)=1$), except for one special vertex. We call a star graph a Type 1
(resp. Type 2) star graph if the color of the special vertex is $1$ (resp. 
$2$). Each Type 1 (resp. Type 2) star graph $\Gamma$ in $G_{d}$ 
is in one to one correspondence with a partition $\sigma_{d}: d=d_{1}+d_{2}+
\cdots +d_{l(\sigma_{d})}$ of an integer 
$d$ whose length $l(\sigma_{d})$ is no less than 2. Obviously, 
we have 
\begin{equation}
|Aut(\Gamma)|=|Aut(\sigma_{d})|. 
\end{equation}
Therefore, the sum of contributions from Type 1 
(resp. Type 2) star graphs in $G_{d}$ is given by,                 
\begin{eqnarray}
&&\frac{1}{d^3}\sum_{\sigma_{d}\in P'_{d}}(-1)^{l(\sigma_{d})+d}
\frac{1}{|Aut(\sigma_{d})|}\prod_{j=1}^{l(\sigma_{d})}
\frac{d\cdot a_{d_{j}}}{j},\;\;\;\qquad\;\;\;(\mbox{Type 1}),\nonumber\\
&&\frac{1}{d^3\cdot f}\sum_{\sigma_{d}\in P'_{d}}(-1)^{l(\sigma_{d})+d}
\frac{1}{|Aut(\sigma_{d})|}\prod_{j=1}^{l(\sigma_{d})}
\frac{d\cdot f\cdot a_{d_{j}}}{j},\;\;\;(\mbox{Type 2}),
\label{star}
\end{eqnarray}
where $P'_{d}$ is the set of partitions $\sigma_{d}$ of the integer $d$ 
whose length $l(\sigma_{d})$ is no less than $2$. 
From elementary combinatorial considerations, the sum of contributions
of all the star graphs can be rewritten by introducing the expression, 
\begin{eqnarray}
s_{d}&:=&\biggl[(-1)^{d}
\exp(-d\bigl(\sum_{j=0}^{\infty}\frac{a_{j}}{j}q^{j}\bigr))
\biggr]_{q^d}+(-1)^{d}a_{d},\nonumber\\
t_{d}&:=&\frac{1}{f}\biggl[(-1)^{d}
\exp(-d\cdot f
\bigl(\sum_{j=0}^{\infty}\frac{a_{j}}{j}q^{j}\bigr))
\biggr]_{q^d}+(-1)^{d}a_{d},
\label{schur} 
\end{eqnarray}
as follows:
\begin{eqnarray}
B(a_1,a_2,a_3,\cdots) &:=&\sum_{d=1}^{\infty}\frac{s_{d}}{d^{3}}q^d, 
\nonumber\\
C(a_1,a_2,a_3,\cdots)&:=&\sum_{d=1}^{\infty}\frac{t_{d}}{d^{3}}q^d. 
\end{eqnarray}
In (\ref{schur}), by $[*]_{q^{d}}$ we mean the operation of picking out the 
coefficient of $q^d$. The notation $B(a_1,a_2,a_3,\cdots)$ and 
$C(a_1,a_2,a_3,\cdots)$ are used to stress that $B$ and $C$ are considered 
as functions in the $a_{d}$'s. It is interesting to note that the expression 
$\biggl[\exp(-d\bigl(\sum_{j=0}^{\infty}\frac{a_{j}}{j}q^{j}\bigr))
\biggr]_{q^d}$ in (\ref{schur}) appears in the process of inverting    
the  power series $w=x+\sum_{d=1}^{\infty}\frac{a_{d}}{d}e^{dx}$, as was 
suggested in \cite{J2}. In other words, summing up graphs with one edge and 
star graphs is closely related to taking the inverse of the mirror map 
in the mirror computation.

Next, we have to sum up the graphs which are neither graphs with one edge
nor star graphs. One can easily see that these graphs are uniquely 
decomposed into the connected sum of star graphs. Then it is a good exercise of 
field theory and combinatorial theory to represent this operation by   
the introduction of a propagator:
\begin{equation}
p_{d}:=(-1)^{d+1}\frac{a_{d}d^{3}}{q^{d}}. 
\end{equation}
We then find the critical values of the action:
\begin{eqnarray}
S:=-\sum_{d=1}^{\infty}\frac{x_{d}y_{d}}{p_{d}}
&+&B(a_1+x_1,a_2+x_2,a_3+x_3,\cdots)-B(a_1,a_2,a_3,\cdots)\nonumber\\
&+&C(a_1+y_1,a_2+y_2,a_3+y_3,\cdots)-C(a_1,a_2,a_3,\cdots),
\end{eqnarray}
where $x_{d}$ and $y_{d}$ are dynamical variables. What remains to be done 
is solving the equation of motion for the $x_{d}$'s and $y_{d}$'s:
\begin{eqnarray}
y_{d}&=&p_{d}\cdot\frac{\partial}{\partial x_{d}}B(a_1+x_1,a_2+x_2,\cdots),
\nonumber\\
x_{d}&=&p_{d}\cdot\frac{\partial}{\partial y_{d}}C(a_1+y_1,a_2+y_2,\cdots).
\label{em}
\end{eqnarray}
At first sight, solving (\ref{em}) appears difficult, but if 
we define recursive formulas:
\begin{eqnarray}
y_{d.n+1}&=&p_{d}\cdot\frac{\partial}{\partial x_{d}}B(a_1+x_{1,n},a_2+x_{2,n},\cdots),
\nonumber\\
x_{d,n+1}&=&p_{d}\cdot\frac{\partial}{\partial y_{d}}C(a_1+y_{1,n},a_2+y_{2,n},\cdots),
\label{em}
\end{eqnarray}
with the initial condition $x_{d,0}=y_{d,0}=0$,  
one can easily obtain the solutions  $x_{d}(a_{*})$ and $y_{d}(a_{*})$ 
in the limit: $\lim_{n\rightarrow\infty}x_{d,n}$ and 
$\lim_{n\rightarrow\infty}y_{d,n}$ respectively. Therefore, by adding up the 
previous contributions from graphs with one edge and star graphs, we finally 
obtain,  
\begin{eqnarray}
F(q,z)&=&A(a_{1},a_{2},a_{3},\cdots)
-\sum_{d=1}^{\infty}\frac{x_{d}(a_{*})y_{d}(a_{*})}{p_{d}}
+B(a_1+x_1(a_{*}),a_2+x_2(a_{*}),a_3+x_3(a_{*}),\cdots)\nonumber\\
&&+C(a_1+y_1(a_{*}),a_2+y_2(a_{*}),a_3+y_3(a_{*}),\cdots).
\end{eqnarray}
Using the above formula, we computed $F(q,1)$ up to degree $10$ and found 
that it is given by,
\begin{equation}
F(q,1):=q + { \frac {1}{8}} q^{2} + { \frac {
1}{27}} q^{3} + { \frac {1}{64}} q^{4} + 
{ \frac {1}{125}} q^{5} + { \frac {1
}{216}} q^{6} + { \frac {1}{343}} q^{7} + 
{ \frac {1}{512}} q^{8} + { \frac {1
}{729}} q^{9} + { \frac {1}{1000}} q^{10}+\cdots, 
\end{equation}
for any $k\geq 1$. This is of course nothing but the multiple cover formula and 
agrees with the results in the previous subsection.
We also write down here the results of $F(q,-1)$ for $k=1,2$ cases:
\begin{eqnarray}
F(q,-1)&=&q - { \frac {7}{8}} q^{2} + { \frac {
55}{27}} q^{3} - { \frac {455}{64}} q^{4} + 
{ \frac {3876}{125}} q^{5} - { 
\frac {33649}{216}} q^{6} + { \frac {296010}{343}
} q^{7} - { \frac {2629575}{512}} q^{8}\nonumber\\ 
&&+ { \frac {23535820}{729}} q^{9} 
 - { \frac {52978783}{250}} q^{10}+\cdots,\;(k=1),\nonumber\\
F(q,-1)&:=&q + { \frac {17}{8}} q^{2} + { 
\frac {325}{27}} q^{3} + { \frac {6545}{64}} q^{
4} + { \frac {135751}{125}} q^{5} + 
{ \frac {2869685}{216}} q^{6} + { 
\frac {61474519}{343}} q^{7} \nonumber\\
&& + { \frac {1329890705
}{512}} q^{8} + { \frac {28987537150}{729}} q^{9} + 
{ \frac {635627275767}{1000}} q^{10}+\cdots,\;(k=2). 
\end{eqnarray}
If we consider $4q\frac{d}{dq}F(q,-1)$, we can see that the result of 
the $k=1$ case completely agrees with the one of the previous 
subsection.
\section{Quantum cohomology of $\p(\mathcal O\oplus \mathcal O (k) \oplus \mathcal O (-2-k))$, $k\ge 1$}

In this section, we detail two methods for determining the isomorphism type of the quantum cohomology ring of $G_k=\p(\mathcal O\oplus \mathcal(k) \oplus \mathcal (-2-k))$ from mirror symmetry. Before going into the particulars, we offer some motivation as to why we are interested in this computation in the first place.

\subsection{Motivation, and a conjecture}

Our initial line of inquiry was the same as above: `How do we describe mirror symmetry for the total space $X_k=\mathcal O(k)\oplus \mathcal O(-2-k)\rightarrow \p^1$?' Before arriving at the equivariant formalism used in Section 3, our first effort was to consider projective bundles, which we now describe. 

In our previous paper \cite{FJ2} we attempted to resolve general questions of local mirror symmetry by replacing noncompact threefolds by projective bundles: 
\begin{eqnarray}
	G_k=\p(\mathcal O\oplus \mathcal O(k) \oplus \mathcal O(-2-k))\rightarrow \p^1.
\end{eqnarray}
In \cite{FJ2}, we then considered the canonical bundle over $G_k$ in order to derive the prepotential. In fact, for spaces with one K\"ahler parameter such as $X_k$, is is actually easier to directly use mirror symmetry for $G_k$. This is the approach we follow here. 

The examples $G_{-1}$ and $G_0$ were considered some time ago by Givental \cite{G2}. We briefly recall his argument. First, we realize $G_{-1}$ as a symplectic quotient $(\C^5-Z)/T^2$, where the weights of the torus action are described by a matrix 
\begin{eqnarray}
	\begin{pmatrix}
	1&1&-1&-1&0 \\
	0&0&1&1&1
	\end{pmatrix}
\end{eqnarray}
and the disallowed locus $Z=\{z_1=z_2=0\} \cup \{z_3=z_4=z_5=0\}$.
With this matrix in hand, we find that the Givental $I$ function for $G_{-1}$ is:
\begin{eqnarray}
I_{-1}=e^{(p_1\log q_1+p_2\log q_2)/\hbar}\sum_{d_1,d_2=0}^{\infty}\frac{q_1^{d_1}q_2^{d_2}\prod_{m=-\infty}^0(-p_1+p_2+m\hbar)^2}{\prod_{m=1}^{d_1}(p_1+m\hbar)^2\prod_{m=1}^{d_2}(p_2+m\hbar)\prod_{m=-\infty}^{-d_1+d_2}(-p_1+p_2+m\hbar)^2}.
\end{eqnarray}
Here, the coefficients of $I_{-1}$ take values in the cohomology ring of $G_{-1}$:
\begin{eqnarray}
	H^*(G_{-1},\C)=\frac{\C[p_1,p_2]}{\langle p_1^2,(-p_1+p_2)^2p_2 \rangle}.
\end{eqnarray}
By expanding this series in powers of $1/\hbar$, we see immediately the Gromov-Witten information of $G_{-1}$:
\begin{eqnarray}
\label{basicI}
I_{-1}=e^{(p_1\log q_1+p_2\log q_2)/\hbar}\Big(1+\frac{Li_2(q_1)p_2^2-2Li_2(q_1)p_1p_2}{\hbar^2}+\frac{q_1q_2+q_2-Li_3(q_1)p_2^3}{\hbar^3}+\dots \Big).
\end{eqnarray}
We note that $I_{-1}$ contains essentially the same information as we find from the Gromov-Witten calculation on the noncompact space $X_{-1}$: a trivial mirror map (i.e. the coefficient of $1/\hbar$ is zero) and the trilogarithm function in $q_1$.

Now let $I_0$ be the $I$ function for $G_0$. Givental shows that (up to a coordinate change by the mirror map) $I_{-1}=I_0$, and hence that $QH^*(G_{-1},\C) \cong QH^*(G_0,\C)$, where $QH^*(X,\C)$ denotes the small quantum cohomology ring of $X$. This implies that $G_0$ has the same Gromov-Witten invariants as $G_{-1}$, which shows that these compactified spaces $G_k$ are in some sense reproducing the physical correspondence between the noncompact $X_{-1}$ and $X_0$. One might therefore hope that this phenomenon would continue to hold for all $G_k$.

In fact, this is true for $k=1,2$, from which we obtained Conjecture \ref{conjecture2}:
For all $k\in \Z$, 
\begin{equation}
QH^*(\p(\mathcal O \oplus \mathcal O(k) \oplus \mathcal O(-2-k))\cong QH^*(\p(\mathcal O \oplus \mathcal O(-1) \oplus \mathcal O(-1)).
\end{equation}       
We will demonstrate this correspondence using only the $I$ function, in the same vein as the above calculation. 

\subsection{Verification of the conjecture using $J$ functions ($k=1,2$)}

First, we need to understand the difference between $G_k$ for $k=-1,0$ and $G_k$ for all other $k\in \Z$. We can exhibit this easily; let $C_1,C_2$ denote the equivalence classes of the base and fiber curve, respectively, and let $p_1, p_2$ be the corresponding K\"ahler classes satisfying $\int_{C_i}p_j=\delta_{ij}$. Then it is easy to see that if $k=-1,0$, we have $c_1(G_k)=3p_2 > 0$. Thus these spaces satisfy the condition of semi-positivity $c_1(X)\ge 0$ used in \cite{G2}.

Now let $k>0$. Recall that we have an equivalence of toric varieties 
\begin{eqnarray}
	\p(\mathcal O\oplus \mathcal O(k) \oplus \mathcal O(-2-k))=\p(\mathcal O\oplus \mathcal O(-k) \oplus \mathcal O(-2-2k)).
\end{eqnarray}
Then we can represent $G_k, k>0$ as a quotient $(\C^5-Z)/T^2$ where the torus action is given by 
\begin{eqnarray}
\label{matrix1}
\begin{pmatrix}
l^1_k \\
l^2_k
\end{pmatrix}=
	\begin{pmatrix}
	1&1&-k&-2-2k&0 \\
	0&0&1&1&1
	\end{pmatrix}
\end{eqnarray}
where $Z=\{z_1=z_2=0\} \cup \{z_3=z_4=z_5=0\}$.
From this matrix, we can compute the first Chern class as a sum of column vectors, which gives $c_1(G_k)=-3kp_1+3p_2$. Clearly, this does not satisfy semi-positivity for $k>0$.

As mentioned in Section 2, for semi-positive manifolds, we have $I \in \C[[h^{-1}]]$, but in general $I\in \C[[h,h^{-1}]]$. In particular, this explains why Givental did not consider $G_k$ for $k>0$ in his original paper, since at the time it was not clear how to remove positive powers of $\hbar$ for the comparison of the $I$ functions.

Now, with the results of \cite{CG}, we are in position to prove the equivalence of quantum cohomology rings of $G_k$ in a similar fashion to \cite{G2}. Let $J_k$ be a $J$ function of $G_k$, $k=1,2$, which is computed from $I_k$ via Birkhoff factorization, as explained in Section 2. From Proposition \ref{jfn-qc}, to show $QH^*(G_k)=QH^*(G_{-1})$, we only need to prove that $J_k=I_{-1}$, up to a coordinate change determined by the coefficient of $\hbar^{-1}$ of $J_k$.
 
We now turn to the computation. Let $I_k$ be the $I$ function of $G_k$. From the weights of the torus action, Eqn. (\ref{matrix1}), we have the following formula for $I_k$:
\begin{eqnarray}
I_k(q,\hbar,\hbar^{-1})=	e^{p\log q/h}\sum_{d_1,d_2}C_{d_1,d_2}(p_1,p_2)q_1^{d_1}q_2^{d_2}
\end{eqnarray}
where the coefficients $C_{d_1,d_2}(p_1,p_2)$ are given by 
\begin{eqnarray}
\frac{\prod_{m=-\infty}^0(N_1p_1+p_2+m\hbar)\prod_{m=-\infty}^0(N_2p_1+p_2+m\hbar)}{\prod_{m=-\infty}^{N_1d_1+d_2}(N_1p_1+p_2+m\hbar)\prod_{m=-\infty}^{N_2d_1+d_2}(N_2p_1+p_2+m\hbar)\prod_{m=1}^{d_1}(p_1+m\hbar)^2\prod_{m=1}^{d_2}(p_2+m\hbar)}
\end{eqnarray}
with $N_1=-k,N_2=-2-2k$. Note that these coefficients take values in the cohomology ring 
\begin{eqnarray}
	H^*(G_k,\C)=\frac{\C[p_1,p_2]}{\langle p_1^2, p_2(N_1p_1+p_2)(N_2p_1+p_2) \rangle}.
\end{eqnarray}
We use the basis $\{1,p_1,p_2,p_1p_2,p_2^2,p_1p_2^2\}$ for $H^*(G_k)$.
Now perform the Birkhoff factorization of $I_k$. Then we acquire a function $J_k$ with asymptotical expansion 
\begin{eqnarray}
	J_k(q,\hbar^{-1})=e^{(p_1 \log q_1+p_2\log q_2)/\hbar}\Big(1+\frac{p_1t_1+p_2t_2}{\hbar}+\frac{p_2^2W_1+p_1p_2W_2}{\hbar^2}+\frac{T_1+p_1p_2^2T_2}{\hbar^3}+\dots\Big).
\end{eqnarray}
Here, $t_1,t_2$ give the mirror map, and the functions $W_1,W_2$ contain the information of Gromov-Witten invariants of this space. Note that this does not exhibit the unusual mirror map behavior mentioned in Section 2; hence, we can proceed directly without need to modify $J_k$. We will see below, for the $F_3$ example, how to deal with the general case.   

 At this point, by simply inserting the inverse mirror map into $W_1$ or $W_2$, we find immediately the usual multiple cover formula for curves.

To complete our present computation, all we have to do is compare the above $J_k$ to the expression for $I_{-1}$ given in Eqn. (\ref{basicI}). This can be done in two steps: (1) Insert the inverse mirror map $q_i(t_j)$ into $J_k$; (2) Make the linear change of cohomology generators 
\begin{equation}
\label{cohshift}
p_1=\tilde p_1, p_2=k\tilde p_1+\tilde p_2.
\end{equation}
 This second step is necessary in order to assure that the cohomology rings for $I_{-1}$ and $J_k$ coincide. After so doing, we arrive at the following expression for $J_k$:
\begin{eqnarray}
e^{(\tilde p_1 \log y_1+ \tilde p_2\log y_2)/\hbar}\Big(1+\frac{\tilde p_2^2 Li_2(y_1y_2^k)-2 \tilde p_1 \tilde p_2 Li_2(y_1y_2^k)}{\hbar^2}+\frac{y_2+y_1y_2^{k+1}-\tilde p_1p_2^2 Li_3(y_1y_2^k)}{\hbar^3}+\dots\Big).
\end{eqnarray}
Here, we have taken $y_i=e^{t_i}$. Note that this is exactly the expected answer: from Eqn. (\ref{matrix1}), $l^1_k+kl^2_k=\begin{pmatrix} 1&1&0&-2-k&k \end{pmatrix}$ is the Calabi-Yau direction, which is why we have the appearance of $y_1y_2^k$ in the above polylog functions.

This completes our proof of the isomorphism $QH^*(G_k)=QH^*(G_{-1})$. Next, we consider the construction of quantum cohomology using connection matrices.
   
\subsection{Alternative proof by connection matrices}
\label{connectionsection}

We can give an alternative proof of the fact that $QH^*(G_k)\cong QH^*(G_{-1})$ by constructing connection matrices for $G_k$, as in \cite{Gu},\cite{I}. Since these  matrices correspond to multiplication in the small quantum cohomology ring, we need only show that the connection matrices on $QH^*(G_k)$ and $QH^*(G_{-1})$ are the same. As in our previous proof, we require only the information of the $I$ function for the computation.

From here forward, we specialize to the case of $G_1=\p(\mathcal O \oplus \mathcal O(-1) \oplus \mathcal O(-4))$. All other cases work out similarly. Then we first carry out the procedure described in Section 2 to compute connection matrices. For this, we start with the fundamental solution 
\begin{eqnarray}
	S^t=\begin{pmatrix} I_1 & \partial_1 I_1 & \partial_2 I_1 & \partial_1\partial_2  I_1 & \partial_2^2 I_1 & \partial_1\partial_2^2 I_1 
	\end{pmatrix}.
\end{eqnarray}
Here $\partial_i=\hbar q_i\partial/\partial q_i$. Then the connection matrices $\Omega_i$ are defined by 
\begin{eqnarray} 
\partial_i S^t=
\begin{pmatrix}\partial_i I_1 & \partial_i \partial_1 I_1 & \partial_i \partial_2 I_1 & \partial_i \partial_1\partial_2  I_1 & \partial_i \partial_2^2 I_1 & \partial_i \partial_1\partial_2^2 I_1 
	\end{pmatrix}^t=\Omega_i	S^t.
\end{eqnarray}
Then we Birkhoff factorize the fundamental solution $S$ as in Section 2, and use the positive part $Q(\hbar)$ to gauge transform the $\Omega_i$. The result is that
the $\hbar$ independent connection matrix $\hat \Omega_1$ corresponding to quantum multiplication by $p_1$ is (up to order 4):
\begin{eqnarray}
\no
	\begin{pmatrix}
	0&1+24q_1q_2+1248q_1^2q_2^2&-4q_1q_2-176q_1^2q_2^2&0&0&0 \\
	0&0&0&-8q_1q_2-340q_1^2q_2^2&q_1q_2+41q_1^2q_2^2&0 \\
	0&0&0&1+20q_1q_2+1084q_1^2q_2^2&-3q_1q_2-135q_1^2q_2^2&0\\
	q_1q_2^2&0&0&0&0& f_1 \\
	-q_1q_2^2&0&0&0&0& f_2 \\
	0&-30q_1q_2^2&5q_1q_2^2&0&0&0
	\end{pmatrix}
\end{eqnarray}
 where $f_1=-4q_1q_2-176q_1^2q_2^2$,$f_2=1+4q_1q_2+368q_1^2q_2^2$.
We also have the matrix $\hat \Omega_2$: 
\begin{eqnarray}
\no
	\begin{pmatrix}
	0&24q_1q_2+1248q_1^2q_2^2&1-4q_1q_2-176q_1^2q_2^2&0&0&0 \\
	0&0&0&1-8q_1q_2-340q_1^2q_2^2&q_1q_2+41q_1^2q_2^2&0 \\
	0&0&0&20q_1q_2+1084q_1^2q_2^2&1-3q_1q_2-135q_1^2q_2^2&0\\
	2q_1q_2^2&0&0&0&0& g_1 \\
	q_2-2q_1q_2^2&0&0&0&0& g_2 \\
	0&5q_2-60q_1q_2^2&10q_1q_2^2&0&0&0
	\end{pmatrix}
\end{eqnarray}
 with $g_1=1-4q_1q_2-176q_1^2q_2^2$, $g_2=5+4q_1q_2+368q_1^2q_2^2$. 
Above, we are using the basis 
\begin{eqnarray}
 \{1,p_1,p_2,p_1p_2,p_2^2,p_1p_2^2\}
 \end{eqnarray}
for $H^*(G_1)$. 

Then in order to show that $QH^*(G_1)\cong QH^*(G_{-1}),$ we need to compare the above matrices $\hat \Omega_1, \hat \Omega_2$ to the connection matrices of $QH^*(G_{-1})$. Let $t_1, t_2$ be the mirror map for $G_1$, as defined by the coefficient of $\hbar^{-1}$ of the function $J_1$ from the previous subsection. Let $q_i(t_j)$ be the inverse mirror map. We perform a coordinate change of the connection matrices $\hat \Omega_i$ via the mirror map:
\begin{eqnarray}
	\tilde \Omega_j=\sum_{i=1}^2 \frac{\partial q_i}{\partial t_j} \hat \Omega_i |_{q=q(t)}
\end{eqnarray}
Finally, we change basis of the $\tilde \Omega_j$ via the linear transformation of Eqn. (\ref{cohshift}). The result of these manipulations are the matrices 

\begin{eqnarray}
	\tilde \Omega_1&=&
	\begin{pmatrix}
	 0&1&0&0&0&0 \\
	 0&0&0&\frac{-2y_{1}y_{2}}{1-y_{1}y_{2}}&
\frac{y_{1}y_{2}}{1-y_{1}y_{2}}&0 \\
	 0&0&0&1&0&0 \\
	 y_1y_2^2&0&0&0&0&0 \\
	 y_1y_2^2&0&0&0&0&1 \\
	 0&-2y_1y_2^2&2y_1y_2^2&0&0&0
	 \end{pmatrix}, \\
	 \tilde \Omega_2&=&
	 \begin{pmatrix}
	 0&0&1&0&0&0 \\
	 0&0&0&1&0&0 \\
	 0&0&0&0&1&0 \\
	 y_1y_2^2&0&0&0&0&1 \\
	 y_2+y_1y_2^2&0&0&0&0&1 \\
	 0&2y_2-2y_1y_2^2&2y_1y_2^2&0&0&0
	 \end{pmatrix}.
\end{eqnarray}
As before, we have set $y_i=e^{t_i}$. One can readily check that the $\tilde \Omega_i$ are the same as those for multiplication in $QH^*(G_{-1})$, thus completing our second proof of the equivalence of quantum cohomology of $G_k$ and $G_{-1}$.

\section{$F_n$ and $K_{F_n}$, $n\ge 3$}

A problem closely related to that of the previous section is the quantum cohomology of $K_{F_n}$, the canonical bundle over the $n$th Hirzebruch surface $F_n=\p(\mathcal O\oplus \mathcal O(-n))$. These spaces have several new features, most notably the presence of a four parameter mirror map for odd $n \ge 3$. In the course of our computations, we arrived at  Conjecture \ref{conjecture3}:
\bigskip

\it
There are two isomorphism types of $QH^*(F_n)$, depending on whether $n$ is odd or even.
\normalfont
\bigskip

In the first subsection, we compute Gromov-Witten invariants and connection matrices for $K_{F_3}$ in a parallel fashion to the previous section, with an emphasis on aspects differing from the $G_k$ examples. The second subsection contains a method for constructing connection matrices for multiplication in the big quantum cohomology ring. Finally, the third subsection is a discussion of the relations which determine quantum cohomology for $F_4$, or equivalently, the differential operators which annihilate $J_{F_4}$.

\subsection{$K_{F_3}$}

A theorem of Coates-Givental \cite{CG} expresses the relationship between the $J$ functions $J_{F_3}$ and $J_{K_{F_3}}$. We therefore begin by computing $J_{F_3}$, as we did in Section 2.

Recall that $F_3$ is the quotient $(\C^4-Z)/T^2$ where the torus action is given by 
\begin{eqnarray}
	\begin{pmatrix}
	1&1&-3&0 \\
	0&0&1&1
	\end{pmatrix}
\end{eqnarray}
and the Stanley-Reisner ideal $Z=\{z_1=z_2=0\} \cup \{z_3=z_4=0\}$.
Then we associate to $F_3$ the $I$ function 
\begin{eqnarray}
	I_{F_3}=e^{(p_1\log q_1 + p_2\log q_2)/\hbar}\sum_d\frac{q_1^{d_1}q_2^{d_2}\prod_{m=-\infty}^0(-3p_1+p_2+m\hbar)}{\prod_{m=-\infty}^{-3d_1+d_2}(-3p_1+p_2+m\hbar)\prod_{m=1}^{d_1}(p_1+m\hbar)^2\prod_{m=1}^{d_2}(p_2+m\hbar)}
\end{eqnarray}
where the coefficients take values in the cohomology of $F_3$: 
\begin{eqnarray}
	H^*(F_3,\C)=\frac{\C[p_1,p_2]}{\langle p_1^2, (-3p_1+p_2)p_2\rangle}.
\end{eqnarray}
As $c_1(F_3)=-p_1+2p_2$, we see that $F_3$ is not semi-positive, and therefore $I_{F_3}\in \C[[\hbar,\hbar^{-1}]]$.

Thus, we must construct the $J$ function as in Definition \ref{birkhofffactor}. The computation leads to the following asymptotic expansion:
\begin{eqnarray}
	J_{F_3}=e^{(p_1\log q_1 + p_2\log q_2)/\hbar}\Big(1+\frac{t_0+t_1p_1+t_2p_2+t_3p_1p_2}{\hbar}+\dots\Big)
\end{eqnarray}
where the $t_i$ are the series 
\begin{eqnarray} \no
	t_0&=&-2q_1q_2-\frac{345}{2}q_1^3q_2-\frac{155209}{3}q_1^5q_2^3-\dots \\ \no
	t_1&=&\frac{135}{2}q_1^2q_2+\frac{181715}{12}q_1^4q_2+\frac{18106223}{3}q_1^6 q_2^3\dots \\ \no
	t_2&=&-16q_1^2q_2-\frac{19267}{6}q_1^4q_2^2-\frac{3619741}{3}q_1^6q_2^3\dots \\
	t_3&=&5q_1+\frac{1901}{3}q_1^3q_2+\frac{2537111}{12}q_1^5q_2^2+\dots
\label{mirror}
	\end{eqnarray}
	Note that here we are using the basis $\{1,p_1,p_2,p_1p_2\}$ for $H^*(F_3)$.
	
	This form of the $J$ function is problematic, for the following reason. In order to determine the mirror map, one looks at the coefficient of $1/\hbar$ of $J_{F_3}$. In the present case, this means that there are 4 power series that determine the mirror map. However, $J_{F_3}$ is only a function of two variables $q_1,q_2$, which means that we must somehow introduce an extra two variables into $J_{F_3}$ in order to invert the power series $t_0\dots t_3$.
	
In what follows, we use the reasoning outlined in Section \ref{nonnef}. Let $\hat \Omega_1, \hat \Omega_2$ be $\hbar$ independent connection matrices which correspond to multiplication by $p_1,p_2$ respectively in $QH^*(F_3)$. Then clearly the identity matrix and $\hat \Omega_1 \hat \Omega_2$ correspond to multiplication by $1$ and $p_1p_2$. Then from Section \ref{nonnef}, we introduce the modified $J$ function
\begin{eqnarray}
	\hat J_{F_3}(q_0\dots q_3,\hbar^{-1})=e^{(q_0I+q_3\hat \Omega_1 \hat \Omega_2)/\hbar}J_{F_3}.
\end{eqnarray}
The function $\hat J_{F_3}$ determines the big quantum cohomology of $F_3$, but that will not concern us here. The important point for the present discussion is that we may now take advantage of the presence of the extra variables to invert the mirror map defined by $t_0\dots t_3$. Set $y_i=e^{t_i}$. Then after inserting the inverse mirror map, we find 
\begin{eqnarray}
\label{realJ}
J_{F_3}'=
	\lim_{y_0,y_3\rightarrow 0}\hat J_{F_3}(y_0\dots y_3, \hbar^{-1})=e^{(p_1\log y_1 + p_2\log y_2)/\hbar}\Big(1+\frac{y_2-2y_1y_2p_1+y_1y_2p_2}{\hbar^2}+\dots\Big).
\end{eqnarray}
Then we note in particular that the $I$ function for the first Hirzebruch surface $F_1$ is given by
\begin{eqnarray}
	I_{F_1}=e^{(\tilde p_1\log q_1 + \tilde p_2\log q_2)/\hbar}\Big(1+\frac{q_2-q_1\tilde p_1+q_1\tilde p_2}{\hbar^2}+\dots\Big)
\end{eqnarray}
and that these two functions agree exactly if we make the substitutions 
\begin{eqnarray}
	q_1=y_1y_2, q_2=y_2, \tilde p_1=p_1, \tilde p_2=p_2-p_1.
\label{3to1}
\end{eqnarray}
Thus, we have demonstrated the equivalence of $QH^*(F_1)$ and $QH^*(F_3)$ at the level of $J$ functions.

Now that we have the correct $J$ function for $F_3$, as defined by $J_{F_3}'$ of eqn. (\ref{realJ}), we can easily compute local Gromov-Witten invariants for $K_{F_3}$. First expand $J_{F_3}'$ in a power series 
\begin{eqnarray}
  J_{F_3}'=\sum_d C'_{d_1,d_2}y_1^{d_1}y_2^{d_2}.	
\end{eqnarray}
Then as in \cite{CG}, we obtain $J_{K_{F_3}}$ by twisting this by a factor corresponding to the canonical bundle of $F_3$ ($K_{F_3}=p_1-2p_2$):
\begin{eqnarray}
	J_{K_{F_3}}=\sum_d C'_{d_1,d_2}y_1^{d_1}y_2^{d_2}\frac{\prod_{m=-\infty}^{0}(p_1-2p_2+m\hbar)}{\prod_{m=-\infty}^{d_1-2d_2}(p_1-2p_2+m\hbar)}.
\end{eqnarray}
Then $J_{K_{F_3}}$ has asymptotics 
\begin{eqnarray}
 e^{(p_1\log y_1 + p_2\log y_2)/\hbar}\Big(1+\frac{s_1p_1+s_2p_2}{\hbar}+\frac{Wp_1p_2}{\hbar^2}\Big)
\end{eqnarray}
where $s_1,s_2$ are the mirror map and $W$ is a function that we use to compute Gromov-Witten invariants.
Notice that the asymptotic expansion terminates at the power $1/\hbar^2,$ which is a feature of $I$ functions for Calabi-Yau spaces.

 Let $\mathcal F$ be the prepotential for $K_{F_3}$. Then the Gromov-Witten invariants can be read off by use of the equation 
\begin{eqnarray}
	W(y_1(s),y_2(s))=-\frac{\partial \mathcal F}{\partial s_1 }+2\frac{\partial \mathcal F}{\partial s_2}.
\end{eqnarray}
 Here $y_i(s)$ is the inverse mirror map. At the end of all this, we arrive at the invariants listed in Table 1. 
\begin{table}[t]
\begin{center}
\begin{tabular}{c|cccccccc}
 &$d_2$&0&1&2&3&4&5&6 \\
\hline
$d_1$& & & & & & & & \\
0& &$N_{0,0}$&-2&0&0&0&0&0\\
1& &0&1&3&5&7&9&11\\
2& &0&$N_{2,1}$&0&0&-6&-32&-110\\
3& &0&0&0&0&0&0&27
 \end{tabular}
\end{center}
\caption{Gromov-Witten invariants for $K_{F_3}$.} 
\end{table}
 
We note that these invariants are the same as those for $K_{F_1}$, except that they appear at a different place on the table \cite{CKYZ}. This is consistent with the results of \cite{CKYZ}, in the sense that the Gromov-Witten invariants found there fore $K_{F_0}$ and $K_{F_2}$ are the same up to their location on the table. Finally, the undetermined invariants $N_{0,0}, N_{2,1}$ cannot be calculated from the mirror symmetric methods we are using here.   
 
\subsection{Connection matrices for $F_3$}
In this section, we compute the connection matrices for the big quantum cohomology ring 
of $F_{3}$ by using the recipe in \cite{Gu}, \cite{I}, \cite{J3}. 
First, by using Birkhoff factorization as described in Subsection 2.3, we can 
construct natural B-model connection matrices associated with $p_{1}$ and 
$p_{2}$,    
\begin{equation}
B_{1}(q_{1},q_{2}):=\left(\begin{array}{cccc}
-2{q_{1}}{q_{2}}-\frac {1035}{2}q_{1}^{3}q_{2}^{2}
& 1 + 135{q_{1}}^{2}{q_{2}}&-32{q_{1}}^{2}{q_{2}}
&\frac {5}{3}{q_{1}}+\frac {1901}{3}{q_{1}}^{3}{q_{2}}\\
10q_{1}^{2}q_{2}^{2}&-864q_{1}^{3}q_{2}^{2}-4{q_{1}}{q_{2}}& 
192q_{1}^{3}q_{2}^{2}+{q_{1}}q_{2}
&- \frac {32}{3}{q_{1}}^{2}{q_{2}}\\ 
 - 12q_{1}^{2}q_{2}^{2}& 1277q_{1}^{3}q_{2}^{2}+3{q_{1}}{q_{2}}& 
- 288q_{1}^{3}q_{2}^{2}-{q_{1}}{q_{2}}&  
\frac {1}{3}+ 13{q_{1}}^{2}{q_{2}}\\
 432{q_{1}}^{3}{q_{2}}^{3}+3{q_{1}}{q_{2}}^{2}& 
-126q_{1}^{2}q_{2}^{2}&30q_{1}^{2}q_{2}^{2}&  
-\frac{1035}{2}q_{1}^{3}q_{2}^{2}-2{q_{1}}
{q_{2}} 
\end{array}\right),
\end{equation}
\begin{equation}
B_{2}(q_{1},q_{2}):=\left(\begin{array}{cccc}
 - 345q_{1}^{3}q_{2}^{2}-2{q_{1}}{q_{2}}&  
\frac {135}{2}{q_{1}}^{2}{q_{2}}& 1-16
{q_{1}}^{2}{q_{2}}&\frac {1901}{9}{q_{1}}^{3}{q_{2}}\\
 10{q_{1}}^{2}{q_{2}}^{2}& -576q_{1}^{3}q_{2}^{2}-4{q_{1}}{q_{2}}&
128q_{1}^{3}q_{2}^{2}+{q_{1}}{q_{2}}&
\frac{1}{3}-\frac{16}{3}{q_{1}}^{2}{q_{2}}\\
{q_{2}}-12{q_{1}}^{2}{q_{2}}^{2}&\frac {2554}{3}q_{1}^{3}q_{2}^{2}+ 
3{q_{1}}{q_{2}}&-192q_{1}^{3}q_{2}^{2}-{q_{1}}{q_{2}}&  
1+\frac{13}{2}q_{1}^{2}{q_{2}}\\
 6{q_{1}}{q_{2}}^{2}+ 432{q_{1}}^{3}{q_{2}}^{3}& 
3{q_{2}}-126{q_{1}}^{2}{q_{2}}^{2}
&30{q_{1}}^{2}{q_{2}}^{2}&- 345q_{1}^{3}q_{2}^{2}-2{q_{1}}{q_{2}}
\end{array}\right),
\end{equation}
up to third order in $q_{1}$. 
In the following, we also use the variables $x_{1}:=\log(q_{1}),\;\;
x_{2}:=\log(q_{2})$. Then what remains to do is to change the 
B-model deformation coordinates $x_{0},x_{1},x_{2},x_{3}$ associated with 
$1,p_{1},p_{2},p_{2}^2$ into A-model flat coordinates 
$t_{0},t_{1},t_{2},t_{3}$. Here, we regard 
$x_{3}$ and $t_{3}$ as the coordinates associated with $p_{2}^{2}$, 
and these differ from those used in previous subsection by a factor of 
$\frac{1}{3}$. In order to execute this coordinate change, we 
have to construct the B-model connection matrices $B_{0}(q_{1},q_{2})$ 
and $B_{3}(q_{1},q_{2})$ for $x_{0}$ and 
$x_{3}$, but these are simply given as follows:
\begin{eqnarray}
B_{0}(q_{1},q_{2})=I_{4},\;\;B_{3}(q_{1},q_{2})=
\bigl(B_{2}(q_{1},q_{2})\bigr)^{2}.
\label{gm}
\end{eqnarray}
The Jacobian matrix between A-model coordinates and B-model coordinates
can be read off from these connection matrices,  
\begin{eqnarray}
\frac{\partial t_{i}}{\partial x_{j}}&=&(B_{j}(q_{1},q_{2}))_{0i},
\end{eqnarray}
and we find that this result completely agrees with the mirror map of  
of Eqn. (\ref{mirror}):
\begin{eqnarray} \no
	t_0&=&-2q_1q_2-\frac{345}{2}q_1^3q_2+\cdots \\ \no
	t_1&=&x_{1}+\frac{135}{2}q_1^2q_2+\cdots \\ \no
	t_2&=&x_{2}-16q_1^2q_2+\cdots \\
	t_3&=&\frac{1}{3}\bigl(5q_1+\frac{1901}{3}q_1^3q_2+\cdots\bigr),
\label{mir2}
	\end{eqnarray}
Therefore, we use (\ref{mir2}) in what follows. With this data, 
we construct intermediate connection matrices $\overline{C}(t_{1},t_{2})$
by the formula: 
\begin{eqnarray}
\overline{C}_{i}(Q_{1},Q_{2})&:=&
\sum_{j=0}^{3}\frac{\partial x_{j}}{\partial t_{i}}B_{j}(q_{1},q_{2}),
\end{eqnarray}
where we have introduced variables $Q_{1}:=\exp(t_{1})$
and $Q_{2}:=\exp(t_{2})$.
The results for $\overline{C}_{1}(Q_{1},Q_{2})$ and 
$\overline{C}_{2}(Q_{1},Q_{2})$ are given by, 
\begin{equation}
\overline{C}_{1}(Q_{1},Q_{2})=\left(\begin{array}{cccc}
0&1&0&0\\
5{Q_1}^{2}{Q_2}^{2}&
-2{Q_1}{Q_2}-\frac{25}{2}{Q_1}^{3}{Q_2}^{2}
&{Q_1}{Q_2}+\frac{25}{2}{Q_1}^{3}{Q_2}^{2}&0\\
10{Q_1}^{2}{Q_2}^{2}&- 2{Q_1}{Q_2}-25{Q_1}^{3}{Q_2}^{2}&
{Q_1}{Q_2}+ 25{Q_1}^{3}{Q_2}^{2}& 
\frac {1}{3}\\ 
3{Q_1}{Q_2}^{2}+\frac {75}{2}{Q_1}^{3}{Q_2}^{3}&
-15{Q_1}^{2}{Q_2}^{2}&15{Q_1}^{2}{Q_2}^{2}&0
\end{array}\right), 
\end{equation}
\begin{equation}
\overline{C}_{2}(Q_{1},Q_{2})=\left(\begin{array}{cccc}
0&0&1&0\\
10{Q_1}^{2}{Q_2}^{2}&
-2{Q_1}{Q_2}-25{Q_1}^{3}{Q_2}^{2}&{Q_1}{Q_2}+ 
25{Q_1}^{3}{Q_2}^{2}&\frac{1}{3}\\
{Q_2}+ 20{Q_1}^{2}{Q_2}^{2}& 
3{Q_1}{Q_2}+\frac {1477}{6}{Q_1}^{3}{Q_2}^{2}& 
{Q_1}{Q_2}+50{Q_1}^{3}{Q_2}^{2}&1\\ 
 6{Q_1}{Q_2}^{2}+\frac {225}{2}{Q_1}^{3}{Q_2}^{3}& 
3{Q_2}-30{Q_1}^{2}{Q_2}^{2}&
30{Q_1}^{2}{Q_2}^{2}&0 
\end{array}\right), 
\end{equation}
where we wrote down the results up to third order in $Q_{1}$.
 At this stage, we have to consider the last 
line of (\ref{mir2}):
\begin{equation}
t_{3}=\frac{1}{3}\bigl(5q_1+\frac{1901}{3}q_1^3q_2+\cdots\bigr)
=\frac{5}{3}Q_{1}+\frac{1777}{18}Q_{1}^{3}Q_{2}+\cdots.
\end{equation}
This means that the B-model expansion point $x_{3}=0$ corresponds 
to $t_{3}=\frac{5}{3}Q_{1}+\frac{1777}{18}Q_{1}^{3}Q_{2}+\cdots$ in 
the A-model coordinate $t_{3}$. Therefore, we have to carry out the parallel 
transport of the $x_{3}$ coordinate by 
$-(\frac{5}{3}Q_{1}+\frac{1777}{18}Q_{1}^{3}Q_{2}+\cdots)$.   
To do this, we have to perturb $\overline{C}_{i}(t_{1},t_{2})\;\;(i=1,2,3)$ by the 
$t_{3}$ coordinate. To this end, we introduce the generating function of 
intermediate Gromov-Witten invariants $w(({\cal O}_{p_{1}})^{n_{1}}
({\cal O}_{p_{2}})^{n_{2}}({\cal O}_{p_{2}^2})^{n_{3}})_{(d_{1},d_{2})}$
as follows: 
\begin{eqnarray}
&&\overline{F}(Q_{1},Q_{2},u_{1},u_{2},u_{3}):=
\frac{1}{2}\sum_{i,j=1}^{3}\eta_{ij}u_{0}u_{i}u_{j}+
\frac{3}{2}(u_{0})^{2}u_{3}\no\\
&&+\sum_{d_{1},d_{2}\geq 0}\sum_{n_{1},n_{2},n_{3}}w(({\cal O}_{p_{1}})^{n_{1}}
({\cal O}_{p_{2}})^{n_{2}}({\cal O}_{p_{2}^2})^{n_{2}})_{(d_{1},d_{2})}
Q_{1}^{d_{1}}Q_{2}^{d_{2}}\frac{u_{1}^{n_{1}}}{n_{1}!}
\frac{u_{2}^{n_{2}}}{n_{2}!}\frac{u_{3}^{n_{3}}}{n_{3}!},
\label{free}
\end{eqnarray}
where $\eta_{ij}$ is the intersection matrix of $F_{3}$:
\begin{equation}
\eta_{ij}  :=  \left(
{\begin{array}{rrrr}
0 & 0 & 0 & 3 \\
0 & 0 & 1 & 0 \\
0 & 1 & 3 & 0 \\
3 & 0 & 0 & 0
\end{array}}
 \right). 
\end{equation}
In (\ref{free}), $n_{1},n_{2},n_{3}$ must satisfy $n_{1}+n_{2}+n_{3} \geq 3$
and $n_{3}=-1-d_{1}+2d_{2}$. The second condition comes from the topological 
selection rule. This function is related to $\overline{C}_{i}(Q_{1},Q_{2})$
by     
\begin{eqnarray}
\partial_{u_{i}}\partial_{u_{j}}\partial_{u_{k}}
\overline{F}(Q_{1},Q_{2},0,0,0)=(\overline{C}_{i}(Q_{1},Q_{2}))_{j}{}^{l}\eta_{lk}.
\end{eqnarray}
As was suggested in \cite{J3}, 
\cite{I},  
$\overline{F}(Q_{1},Q_{2},u_{1},u_{2},u_{3})$ can be fully determined by the modified 
K\"ahler equations:
\begin{eqnarray}
&&\frac{\partial}{\partial u_{1}}\overline{F}(Q_{1},Q_{2},u_{1},u_{2},u_{3})=
\frac{\partial}{\partial t_{1}}\overline{F}(Q_{1},Q_{2},u_{1},u_{2},u_{3})-
\frac{\partial t_{3}}{\partial t_{1}}\frac{\partial}{\partial u_{3}}
\overline{F}(Q_{1},Q_{2},u_{1},u_{2},u_{3}),
\nonumber\\
&&\frac{\partial}{\partial u_{2}}\overline{F}(Q_{1},Q_{2},u_{1},u_{2},u_{3})=
\frac{\partial}{\partial t_{2}}\overline{F}(Q_{1},Q_{2},u_{1},u_{2},u_{3})-
\frac{\partial t_{3}}{\partial t_{2}}\frac{\partial}{\partial u_{3}}  
\overline{F}(Q_{1},Q_{2},u_{1},u_{2},u_{3}),\nonumber\\
\end{eqnarray}
and the associativity equation:
\begin{eqnarray}
&&\partial_{u_{i}}\partial_{u_{j}}\partial_{u_{k}}
\overline{F}(Q_{1},Q_{2},0,0,u_{3})\eta^{kl}
\partial_{u_{l}}\partial_{u_{m}}\partial_{u_{n}}\overline{F}
(Q_{1},Q_{2},0,0,u_{3})
\no\\
&&=\partial_{u_{i}}\partial_{u_{m}}\partial_{u_{k}}
\overline{F}(Q_{1},Q_{2},0,0,u_{3})\eta^{kl}
\partial_{u_{l}}\partial_{u_{j}}\partial_{u_{n}}
\overline{F}(Q_{1},Q_{2},0,0,u_{3}),
\end{eqnarray}
where $\eta^{ij}$ is the inverse matrix of $\eta_{ij}$. 
Then the perturbed intermediate connection matrices are given by,
\begin{equation}
(\overline{C}_{i}(Q_{1},Q_{2},u_{1},u_{2},u_{3}))_{j}{}^{l}
=\partial_{u_{i}}\partial_{u_{j}}\partial_{u_{k}}
\overline{F}(Q_{1},Q_{2},u_{1},u_{2},u_{3})\eta^{kl}.
\end{equation}
Finally, we can construct A-model connection matrices for $F_{3}$
by the parallel transport, 
\begin{equation}
C_{i}(Q_{1},Q_{2})=\overline{C}_{i}(Q_{1},Q_{2},0,0,-t_{3}).
\end{equation}
Below, we write down the exact results of  $C_{1}(Q_{1},Q_{2})$
and $C_{2}(Q_{1},Q_{2})$. 
\begin{equation}
C_{1}(Q_{1},Q_{2}):=\left(\begin{array}{cccc}
0 & 1 & 0 & 0 \\
0 &-2{Q_1}{Q_2} &{Q_1}{Q_2} & 0\\
0 &-2{Q_1}{Q_2} &{Q_1}{Q_2} & \frac {1}{3}\\ 
3{Q_1}{Q_2}^{2}& 0 & 0 & 0
\end{array}\right), 
\end{equation} 
\begin{equation}
C_{2}(Q_{1},Q_{2}):=\left(\begin{array}{cccc}
0 & 0 & 1 & 0 \\
0 &-2{Q_1}{Q_2} &{Q_1}{Q_2} &\frac {1}{3}\\ 
{Q_2}&-2{Q_1}{Q_2} &{Q_1}{Q_2} & 1 \\
6{Q_1}{Q_2}^{2} & 3{Q_2} & 0 & 0
\end{array}\right). 
\end{equation} 
By applying the coordinate change (\ref{3to1}) to the above results, 
we can obtain exactly the same connection matrices as those for 
$QH^{*}(F_{1})$. 
We have thus verified our conjecture in the $F_{3}$ case.

\subsection{$F_4$:  quantum differential equations}

We also consider $F_4$, but our focus is a bit different from that of the previous sections. We will show that $QH^*(F_4)=QH^*(F_2)$, but here, this will be done by making use the relations which determine quantum cohomology for $F_4$. As we will see, we cannot simply use the Picard-Fuchs equations as relations on quantum cohomology for toric $X$ not satisfying $c_1(X)\ge 0$. Thus, the basic question we are exploring is: What happens to the Picard-Fuchs equations when we perform Birkhoff factorization? 

As usual, we represent $F_4=(\C^4-Z)/T^2$ with torus action 
\begin{eqnarray}
	\begin{pmatrix}
	1&1&-4&0 \\
	0&0&1&1
	\end{pmatrix}
\end{eqnarray}
with $Z=\{z_1=z_2=0\} \cup \{z_3=z_4=0\}$.
The $I$ function is thus 
\begin{eqnarray}
	I_{F_4}=e^{(p_1\log q_1 + p_2\log q_2)/\hbar}\sum_d\frac{q_1^{d_1}q_2^{d_2}\prod_{m=-\infty}^0(-4p_1+p_2+m\hbar)}{\prod_{m=-\infty}^{-4d_1+d_2}(-4p_1+p_2+m\hbar)\prod_{m=1}^{d_1}(p_1+m\hbar)^2\prod_{m=1}^{d_2}(p_2+m\hbar)}
\end{eqnarray}
where the generators of cohomology lie in
\begin{eqnarray}
	H^*(F_4,\C)=\frac{\C[p_1,p_2]}{\langle p_1^2, (-4p_1+p_2)p_2\rangle}.
\end{eqnarray}

We want to consider the differential equations which annihilate $I_{F_4}$. From \cite{G2}, these are 
\begin{eqnarray} \no
	\mathcal D_1&=&\theta_1^2-q_1(-4\theta_1+\theta_2)(-4\theta_1+\theta_2-\hbar)(-4\theta_1+\theta_2-2\hbar)(-4\theta_1+\theta_2-3\hbar), \\
\mathcal D_2&=&\theta_2(-4\theta_1+\theta_2)-q_2	
\end{eqnarray}
where $\theta_i=\hbar q_i \partial/ \partial q_i$. Now, if $F_4$ was a semi-positive manifold, we would be able to represent quantum cohomology simply by considering 
\begin{eqnarray}
	\frac{\C[\theta_1,\theta_2,q_1,q_2,\hbar]}{\langle \mathcal D_1, \mathcal D_2 \rangle}.
\end{eqnarray}
However, in the present case the Gr\"obner basis calculation is intractable. The reason is essentially that the coefficient of $q_1$ of $\mathcal D_1$ contains order 4 terms, while the ordinary cohomology algebra of $F_4$ is only two dimensional. 

In light of the calculations of the previous sections, the following resolution presents itself. We observe that the differential operator $\mathcal D_1$ contains higher order powers of $\theta_1,\theta_2$ exactly because the asymptotic expansion of $I_{F_4}$ has positive powers of $\hbar$. Then we know from \cite{CG} that we can eliminate positive powers of $\hbar$ simply by Birkhoff factorizing. This suggests that we can find a normalized form of $\mathcal D_1,\mathcal D_2$ by computing instead differential operators $ \hat{\mathcal D_1}, \hat{\mathcal D_2}$ such that
\begin{eqnarray}
	\hat{\mathcal D_i} J_{F_4}=0.
\end{eqnarray}
  
  This can be done easily. First Birkhoff factorize $I_{F_4}$ to obtain $J_{F_4}$, and then insert the inverse mirror map into $J_{F_4}$. Here the mirror map is determined by the coefficient of $1/\hbar$ in the asymptotic expansion of $J_{F_4}$:
\begin{eqnarray}
	J_{F_4}(q_1,q_2,\hbar^{-1})=1+\frac{p_1t_1(q)+p_2t_2(q)}{\hbar}+\dots
\end{eqnarray}
  
We observe that for $F_4$, as for the $\p(\mathcal O\oplus \mathcal O(k)\oplus \mathcal O(-2-k))$ examples, we do not need to introduce extra variables to $J_{F_4}$, since there are only 2 functions $t_1,t_2$ which determine the mirror map. Let $q_i(t)$ denote the inverse mirror map.
  
Then we directly compute the annihilators of $J_{F_4}(q_1(t),q_2(t),\hbar^{-1})$ to be 
\begin{eqnarray}
	\hat{\mathcal D_1}&=&\hat{\theta_1^2}-y_1y_2^2, \\
	\hat{\mathcal D_2}&=&\hat{\theta_2}(\hat{\theta_2}-4\hat{\theta_1})+y_2(4y_1y_2-1).
\end{eqnarray}
Here $y_i=e^{t_i}$ and $\hat{\theta_i}=\hbar y_i\partial/\partial y_i$. These new operators have none of the problems of $\mathcal D_1,\mathcal D_2$ and hence we can write the quantum cohomology ring of $F_4$ as:
\begin{eqnarray}
	QH^*(F_4,\C)=\frac{\C[\hat{\theta_1},\hat{\theta_2},y_1,y_2]}{\langle \hat{\mathcal D_1},\hat{\mathcal D_2}\rangle }.
\end{eqnarray}
 This is exactly as expected, if one compares with the computation for $F_2$ from \cite{Gu}, and furthermore proves that $QH^*(F_4)=QH^*(F_2)$. 

In closing, we mention that we will \it always \normalfont be able to find a  well-behaved set of differential operators annihilating $J$. The reason is simply that $J$ satisfies the relations 
\begin{eqnarray}
	\hbar \frac{\partial}{\partial t_a}\frac{\partial}{\partial t_b}J_X(t,\hbar^{-1})=\sum_c A^c_{ab}(t)\frac{\partial}{\partial t_c}J_X(t,\hbar^{-1}).
\end{eqnarray}

\section{Conclusion}

In this paper, we have developed a complete computational scheme for determining Gromov-Witten invariants and quantum cohomology rings for $X_k=\mathcal O(k)\oplus \mathcal O(-2-k)\rightarrow \p^1, \ k \ge 1,$ as well as non-nef toric varieties, by using mirror symmetry. Several new features have emerged in the course of our study. For $X_k$, we have seen, first of all, that we need to work in an equivariant setting to compute the correct Gromov-Witten invariants. The second new aspect is the necessity of using Birkhoff factorization to deal with the $I$ function for $X_k$. This second point is not immediately obvious, but given the nature of $I$ functions for the spaces $\p(\mathcal O\oplus \mathcal O(k)\oplus \mathcal O(-2-k))$, the introduction of Birkhoff appears more naturally.

Several questions are raised by our work here. The first is the behavior of Picard-Fuchs equations across the Birkhoff factorization. We know how to derive the $J$ function from $I$ via Birkhoff, but the corresponding transformation of Picard-Fuchs equations is less clear. Another question is the meaning of the Gromov-Witten invariants of $X_k$ with the anti-diagonal action. This corresponds to the equivariant Calabi-Yau case, and as such agrees with results from physics, but we are not aware of the actual geometric meaning of these numbers. At any rate, we hope to address these and other questions in future work.

\end{document}